\documentclass[12pt]{article}
\usepackage{amsmath,amsfonts}
\usepackage{mathrsfs}
\usepackage{latexsym} %amsmath,
\newcommand{\nref}[1]{(\ref{#1})}
\def\Cal{\cal}

\newtheorem{remark}{Remark}[section]
\def\text#1{\hbox{#1}}

\def\endproof{\mbox{\ $\Box$}}
\def\Cal{\cal}

\newcommand{\R}        {{{\rm I\! R}}}

\def\1{\mbox{1\hspace{-.25em}I}}

\newcommand{\CA}{{\Cal{A}}}

\newcommand{\CN}{{\Cal{N}}}

\newcommand{\eqdef}{{ \,\stackrel{\Delta}{=}\, }}

\def\e{\varepsilon}
\def\b{\beta}
\def\g{\gamma}

\def\l{\left}
\def\r{\right}

\def\Var{\rm {Var}}

\def\CN{{\cal{N}}}

\def\bX{{\bf X}}
\def\bY{{\bf Y}}

\def\Var{{\rm Var}}

\newtheorem{theorem}{Theorem}[section]

\newtheorem{lemma}{Lemma}[section]

\newtheorem{corollary}{Corollary}[section]
\numberwithin{equation}{section}

\begin{document}

\title{\bf Sparse classification boundaries}
\author{Yuri I. Ingster$^1$, Christophe Pouet$^2$  and Alexandre B. Tsybakov
$^3$}

\maketitle

\begin{abstract}Given a
training sample of size $m$ from a $d$-dimensional population, we
wish to allocate a new observation $Z\in \R^d$ to this population or
to the noise. We suppose that the difference between the
distribution of the population and that of the noise is only in a
shift, which is a sparse vector. For the Gaussian noise, fixed
sample size $m$, and the dimension $d$ that tends to infinity, we
obtain the sharp classification boundary and we propose classifiers
attaining this boundary. We also give extensions of this result to
the case where the sample size $m$ depends on $d$ and satisfies the
condition $(\log m)/\log d \to \gamma$, $0\le \gamma<1$, and to the
case of non-Gaussian noise satisfying the Cram\'er condition.
\end{abstract}


\thanks{$^1$Research partially supported by the RFBI Grant 08-01-00692-a and by Grant
NSh--638.2008.1.}

\thanks{$^2$Research partially supported by the grant ANR-07-BLAN-0234 and
 by PICS--2715.}

\thanks{$^3$Research partially supported by the grant ANR-06-BLAN-0194, by the PASCAL
Network of Excellence and Isaac Newton Institute for Mathematical
Sciences in Cambridge (Statistical Theory and Methods for Complex,
High-Dimensional Data Programme, 2008).}

\section{Introduction}
\subsection{Model and problem}
 Let ${\bX}=(X_1,\dots,X_n)$
and ${\bY}=(Y_1,\dots,Y_m)$ be two i.i.d. samples from two different
populations with probability distributions $P_X$ and $P_Y$ on $\R^d$
respectively. Here
$$
X_i=(X_i^1,\dots,X_i^d),\quad Y_j=(Y_j^1,\dots,Y_j^d)
$$
where $X_i^k$ and $Y_j^k$ are the components of $X_i$ and $Y_j$. We
consider the problem of discriminant analysis when the dimension of
the observations $d$ is very large (tends to $+\infty$). Assume that
we observe a random vector $Z=(Z^1,\dots,Z^d)$ independent of
$({\bX}, {\bY})$ and we know that the distribution of $Z$ is either
$P_X$ or $P_Y$. Our aim is to classify $Z$, i.e., to decide whether
$Z$ comes from the population with distribution $P_X$ or from that
with distribution $P_Y$.

In this paper we assume that
\begin{equation}\label{model1}
X^k_i=v_k + \xi^k_i,\qquad \qquad Y^k_j= u_k+\eta^k_j,
\end{equation}
where $v=(v_1,\dots,v_d),\ u=(u_1,\dots,u_d)$ are deterministic mean
vectors and the errors $\xi_i^1,\dots,\xi_i^d$, $\eta_j^1,\dots,
\eta_j^d$ are (unless other conditions are explicitly mentioned)
jointly i.i.d. zero mean random variables with probability density
$f$ on $\R$.

Distinguishing between $P_X$ and $P_Y$ presents a difficulty only
when the vectors $v$ and $u$ are close to each other. A particular
type of closeness for large $d$ can be characterized by the sparsity
assumption \cite{IS02a,DJ04}
 that we shall adopt in this paper. As in \cite{IS02a,DJ04}, we introduce the following set of sparse
vectors in $\R^d$ characterized by a positive number $a_d$ and a
{\it sparsity index} $\beta\in (0,1]$:
$$
U_{\beta, a_d}= \left\{u=(u_1,\dots,u_d):\ u_k=a_d\e_k ,\quad
\e_k\in \{0,1\},\quad cd^{1-\b} \le \sum_{k=1}^d\e_k \le
Cd^{1-\b}\right\}.
$$
Here $0<c < C<+\infty$ are two constants that are supposed to be
fixed throughout the paper. %; for $\b=1$ we suppose $c\ge 1$.
The value $p=d^{-\beta}$ can be interpreted as the ``probability" of
occurrence of non-zero components in vector $u$.

In what follows we shall deal only with a special case of model
(\ref{model1}) that was also considered recently by \cite{HPG}. Namely, we assume:
$$
v=0, \quad u\in U_{\b,a_d}.
$$
In this paper we establish the classification boundary, i.e., we
specify the necessary and sufficient conditions on $\beta$ and $a_d$
such that successful classification is possible. Let us first define
the notion of successful classification. We shall need some
notation. Let $\psi$ be a decision rule, i.e., a measurable function
of $\bX,\bY,Z$ with values in $[0,1]$. If $\psi=0$ we allocate $Z$
to the $P_X$-population, whereas for $\psi=1$ we allocate $Z$ to the
$P_Y$-population. The rules $\psi$ taking intermediate values in
$(0,1)$ can be interpreted as randomized decision rules. Let
$P_{H_0}^{(u)}$ and $P_{H_1}^{(u)}$ denote the joint probability
distributions of $\bX,\bY,Z$ when $Z\sim P_X$ and $Z\sim P_Y$
respectively, and let $E_{H_0}^{(u)}$, ${\rm Var}_{H_0}^{(u)}$ and
$E_{H_1}^{(u)}$, ${\rm Var}_{H_1}^{(u)}$ denote the corresponding
expectation and variance operators. We shall also denote by
$P^{(u)}$ the distribution of ${\bf Y}$ and by $E^{(u)}$, ${\rm
Var}^{(u)}$ the corresponding expectation and variance operators.
Consider the Bayes risk
$$
{\cal R}_B(\psi)= \pi E_{H_0}^{(u)} (\psi) + (1-\pi)
E_{H_1}^{(u)}(1-\psi),
$$
where $0<\pi<1$ is a prior probability of the $P_X$-population, and
the maximum risk
$$
{\cal R}_M(\psi)= \max \Big(E_{H_0}^{(u)} (\psi),
E_{H_1}^{(u)}(1-\psi)\Big).
$$
Let ${\cal R}(\psi)$ be either the Bayes risk ${\cal R}_B(\psi)$ or
the maximum risk ${\cal R}_M(\psi)$.

We shall say that {\it successful classification is possible} if
$\beta$ and $a_d$ are such that
\begin{equation}\label{bs}
\lim_{d\to +\infty} \inf_{\psi} \sup_{u\in U_{\b,a_d}} {\cal R}(\psi)
= 0
\end{equation}
for ${\cal R}={\cal R}_M$ and ${\cal R}={\cal R}_B$ with any fixed
$0<\pi<1$. Conversely, we say that {\it successful classification is
impossible} if $\beta$ and $a_d$ are such that
\begin{equation}\label{bi}
\liminf_{d\to +\infty} \inf_{\psi} \sup_{u\in U_{\b,a_d}} {\cal
R}(\psi) ={\cal R}_{max},
\end{equation}
where ${\cal R}_{max}=1/2$ for ${\cal R}={\cal R}_M$ and ${\cal
R}_{max}=\min(\pi,1-\pi)$ for ${\cal R}={\cal R}_B$ with $0<\pi<1$.

We call (\ref{bs}) the {\it upper bound of classification} and
(\ref{bi}) the {\it lower bound of classification}. The lower bound
(\ref{bi}) for the maximum risk ${\cal R}={\cal R}_M$ is interpreted
as the fact that no decision rule is better (in a minimax sense)
than the simple random guess. For the Bayes risk ${\cal R}_B$, the
lower bound \nref{bi} is attained at the degenerate decision rule
that does not depend on the observations: $\psi\equiv0$ if $\pi>1/2$ or
$\psi\equiv1$ if $\pi\le 1/2$.
%This is why the quantity ${\cal
%R}_{max}$ appears on the right-hand side of (\ref{bi}).

The condition on $(\beta, a_d)$ corresponding to the passage from
(\ref{bs}) to (\ref{bi}) is called the {\it classification
boundary}. We shall say that a classifier $\psi=\psi_d$ is {\it
asymptotically optimal} (or that $\psi$ {\it attains the
classification boundary}) if, for all $\beta$ and $a_d$ such that
successful classification is possible, we have
%{\bf (I am not sure that it is
%good definition in the point "for all", we need some "uniformity" in
%the conditions, at least for $\b\in (0,1/2]$. I'll try discuss it in
%Sections \ref{sec2} and \ref{sec_ad} (on adaptivity)) },
\begin{equation}\label{bs11}
\lim_{d\to +\infty} \sup_{u\in U_{\b,a_d}} {\cal R}(\psi) = 0
\end{equation}
where ${\cal R}={\cal R}_M$ or ${\cal R}={\cal R}_B$ with any fixed
$0<\pi<1$.

\subsection{Main results}

According to the value of $\beta$, we shall distinguish between {\it
moderately sparse vectors} and {\it highly sparse vectors}. This
division depends on the relation between $m$ and~$d$. For $m$ not
too large, i.e., when $\log m=o(\log d)$, moderately sparse vectors
correspond to $\beta\in(0,1/2]$ and highly sparse vectors to
$\beta\in(1/2, 1)$. For large $m$, i.e., when $\log m\sim \gamma
\log d,\ \g\in (0,1)$, moderate sparsity corresponds to $\b\in
(0,(1-\g)/2]$ and high sparsity to $\b\in ((1-\g)/2, 1-\g)$.

The classification boundary for moderately sparse vectors is
obtained in a relatively simple way (cf. Section 2).   It is of the
form
\begin{equation}\label{bound1}
R_d\eqdef d^{1/2-\b}a_d\asymp 1.
\end{equation}
This means that successful classification is possible if $R_d \to
+\infty$, and it is impossible if $R_d \to 0$ as $d\to +\infty$. The
result is valid  both for $\beta\in(0,1/2]$ and $m\ge 1$ fixed or
for  $m$ depending on  $d$ such that $\log m\sim \gamma \log d,\
\g\in (0,1)$ as $d\to +\infty$ and $\beta\in(0,(1-\g)/2]$. Moreover,
(\ref{bound1}) holds under weak assumptions on the noise. In
particular, for the upper bound of classification we only need to
assume that the noise has mean zero and finite second moment (cf.
Section 2). The lower bound is proved under a mild regularity
condition on the density $f$ of the noise.

The case of highly sparse vectors is more involved. We establish the
classification boundary for the following scenarios:
\begin{itemize}
\item[(A)]  $m\ge 1$ is a fixed integer, and the noise density $f$ is
Gaussian $\CN(0,\sigma^2)$ with known or unknown
 $\sigma>0$;
\item[(B)]  $m\to+\infty$ as $d\to +\infty$, $\log m=o(\log d)$, and $f$ is
Gaussian $\CN(0,\sigma^2)$ with known or unknown $\sigma>0$.
\item[(C)]  $\log m\sim \gamma
\log d,\ \g\in (0,1)$, and $f$ is Gaussian $\CN(0,\sigma^2)$ with
known or unknown $\sigma>0$.
\end{itemize}
The upper bounds are extended to the following additional scenario:
\begin{itemize}
\item[(D)]   $m\to+\infty$ as $d\to +\infty$, $\log m\sim\g \log d,\ 0\le\g<1,\ m/\log d\to +\infty$,
and the noise satisfies the Cram\'er condition.
\end{itemize}

The conditions on the noise in (A)--(D) are crucial and, as we shall
see later, they suggest that a special dependence of $a_d$ on $d$
and $m$ of the form $a_d\asymp \sqrt{(\log d)/m}$ is meaningful in
the highly sparse case. More specifically, we take
\begin{equation}\label{ad}
a_d = s\sigma\sqrt{\log d},\quad x_1= s\sqrt{m+1},
\end{equation}
where $x_1>0$ is fixed. The classification boundary in (A, B, D) is
then expressed by the following condition on $\beta$, $s$ and $m$:
\begin{equation}\label{detb}
x_1 = \phi(\beta)
\end{equation}
where
\begin{equation}\label{phi}
\phi(\b)= \begin{cases}
\phi_1(\b) & \rm{if} \ 1/2 < \b\le 3/4,\\
\phi_2(\b) & \rm{if} \ 3/4 < \b < 1,
\end{cases}
\end{equation}
with
\begin{equation}\label{phi12}
\phi_1(\b)= \sqrt{2\b-1},\quad
\phi_2(\b)=\sqrt{2}\Big(1-\sqrt{1-\b}\Big).
\end{equation}
In other words, successful classification is possible if $x_1 \ge
\phi(\beta)+\delta$, and it is impossible if $x_1 \le
\phi(\beta)-\delta$,  for any $\delta>0$ and $d$ large enough. This
classification boundary is also extended to the case where $x_1$
depends on $d$ but stays bounded.

For Scenario (C) let $a_d=\sigma x\sqrt{(\log d)/m}$ with fixed
$x>0$. We show that in this framework successful classification is
impossible if $\b>1-\g$ (cf. 1$^\circ$ in Section 2), and therefore
we are interested in $\b\in ((1-\g)/2, 1-\g)$. Set
$\b^*=\b/(1-\g)\in (1/2,1)$ and $x^*=x/\sqrt{1-\g}$. Then the
classification boundary is of the form
$$
x^* = \phi(\b^*),
$$
for the function $\phi(\b)$ defined above.

Note that if $f$ is known, the distribution $P_X$ is also known.
This means that we do not need the sample $\bX$ to construct
decision rules. Thus, in Scenarios (A), (B) and (C) when $\sigma$ is
known we can suppose w.l.o.g. that only the sample $\bY$ is
available; this remark remains valid in the case of unknown
$\sigma$, as we shall see it later. As to Scenario (D), we shall
also treat it under the assumption that only the sample $\bY$ is
available (w.l.o.g. if $f$ is known), to be consistent with other
results. However, if $f$ is not known, the sample $\bX$ contains
additional information which can be used. The results for this case
under Scenario (D) are similar to those that we obtain below but
they are left beyond the scope of the paper.

For $m=0$ (i.e., when there is no sample $\bY$) the problem that we
consider here reduces to the problem of signal detection in growing
dimension $d$, cf.  \cite{I,IS01a,IS01b,IS02a,IS02b,DJ04,JW},
 and our classification boundary coincides with the {\it
detection boundary} established in \cite{I}. Sharp asymptotics in
the detection problem was studied in \cite{I} (see also
\cite{IS02a}, Chapter 8)
for known $a_d$ or  $\b$. Adaptive problem
(this corresponds to unknown $a_d$ and $\b$) was studied in
\cite{IS01a,IS01b}. Various procedures attaining the detection boundary
were proposed in \cite{IS02b,DJ04,JW}. Ingster and Suslina \cite{IS02b}
 introduced a method attaining the detection boundary
based on the combination of three different procedures for the zones
$\beta\in(0,1/2]$, $\beta\in(1/2,3/4]$ and $\beta\in(3/4,1)$. Later
Donoho and Jin \cite{DJ04} showed that a test based on the higher
criticism statistic attains the detection boundary simultaneously
for these zones. More recently Jager and Wellner \cite{JW} proved that
the same is true for a large class of statistics including the
higher criticism statistic.

The paper of Hall {\it et al.} \cite{HPG} deals with the same classification model as
the one we consider here but study a problem which is different from
ours. They analyse the conditions under which some simple (for
example, minimum distance) classifiers $\psi$ satisfy
\begin{equation}\label{bshall}
\lim_{d\to +\infty} E_{H_0}^{(u)}(\psi) = 0.
\end{equation}
Hall {\it et al.} \cite{HPG} conclude that for minimum distance
classifiers (\ref{bshall}) holds if and only if $0<\beta<1/2$. This
implies that such classifiers cannot be optimal for $1/2\le
\beta<1$. They also derive (\ref{bshall}) for some other classifiers
in the case $m=1$.

The results of this paper and their extensions to the multi-class
setting were summarized in \cite{Po} and presented at the Meeting
``Rencontres de Statistique Math\'ematique" (Luminy, December 16-21,
2008) and at the Oberwolfach meeting ``Sparse Recovery Problems in
High Dimensions: Statistical Inference and Learning Theory" (March
15-21, 2009). In a work parallel to ours, Donoho and Jin
\cite{DJ08a,DJ08b} and Jin \cite{jin2009} independently and
contemporaneously have analysed a setting less general than the
present one. They did not consider a minimax framework, but rather
demonstrated that the higher criticism (HC) methodology can be
successfully extended to the classification problem.  Donoho and Jin
\cite{DJ08b} showed that, for a special case of Scenario (B), the
``ideal" HC statistic attains the same upper bound of classification
that we prove below. Together with our lower bound, this implies
that the ``ideal" HC statistic is asymptotically optimal, in the
sense defined above, for the Scenario (B). Donoho and Jin announce
that similar results for the HC statistic in Scenarios (A) and (C)
will appear in their work in preparation.

This paper is organized as follows. Section 2 contains some
preliminary remarks. In Section~3 we present the classification
boundary and asymptotically optimal classifier for moderately sparse
vectors under rather general conditions on the noise. In Section 4
we give the classification boundary and asymptotically optimal
classifiers for highly sparse vectors under Scenarios (A), (B) and (C).
%(with known$\sigma$)
Section 5 provides an extension
%to unknown $\sigma$ and
to Scenario (D). Proofs of the lower and upper bounds of
classification are given in Sections 6 and 7, respectively.
%%%%%%%%%%%%%%

\section{Preliminary remarks}\label{pre}

In this section we collect some basic remarks on the problem
assuming that $f$ is the standard Gaussian density. As a starting
point, we discuss some natural limitations for $a_d$.

1$^\circ$. Remark that $a_d$ cannot be too small. Indeed, assume
that instead of the set $U_{\beta,a_d}$ we have only one vector
$u=(a_d\e_1,\dots,a_d\e_d)$ with known $\e_k\in\{0,1\}$. Then we get
a familiar problem of classification with two given Gaussian
populations. The notion of classification boundary can be defined
here in the same terms as above, and the explicit form of the
boundary can be derived from the standard textbook results. It is
expressed through the behavior of $ Q^2_d \eqdef  a_d^2
\sum_{k=1}^d\e_k$:
\begin{itemize}
\item if $Q_d\to 0$, then successful classification is
impossible: $$\liminf_{d \to +\infty} \inf_{\psi} {\cal
R}(\psi)={\cal R}_{max},$$
\item if $Q_d\to+\infty$, then
successful classification is realized by the maximum likelihood
classifier $\psi^*=\1_{\{T^*>0\}}$ where
$$
T^*=\sum_{k=1 : \e_k=1}^d (Z^k-a_d/2).
$$
\end{itemize}
Here and below $\1_{\{\cdot\}}$ denotes the indicator function.

If we assume that $\sum_{k=1}^d\e_k\asymp d^{1-\beta}$, we
immediately obtain some consequences for our model defined in
Section 1. We see that successful classification in that model is
impossible if $a_d$ is so small that $d^{1-\beta}a_d^2=o(1)$, and it
makes sense to consider only such $a_d$ that
\begin{equation}\label{0}
d^{1-\beta}a_d^2\to+\infty.
\end{equation}
In particular, for $\g>1-\b$ successful classification is impossible
under Scenario (C) with $a_d\asymp \sqrt{(\log d)/m}$.

We shall see later that (\ref{0}) is a rough condition, which is
necessary but not sufficient for successful classification in the
model of Section 1. For example, in that model with $\beta\in
(0,1/2]$ and fixed $m$, the value $a_d$ should be substantially
larger than given by the condition $d^{1-\beta}a_d^2\asymp 1$, cf.
(\ref{bound1}).

2$^\circ$. Our second remark is that, on the other extreme, for
sufficiently large $a_d$ the problem is trivial. Specifically,
non-trivial results can be expected only under the condition
\begin{equation}\label{1*}
x\eqdef a_d\sqrt{m/\log d }\le 2\sqrt{2},
\end{equation}
Indeed, assume that
\begin{equation}\label{1}
x> 2\sqrt{2}.
\end{equation}
Then the problem becomes simple in the sense that successful
classification is easily realisable under (\ref{1}) and the
classical condition (\ref{0}). Indeed, take an analog of the
statistic $T^*$ where $a_d$ and $\e_k$ are replaced by their natural
estimators:
\begin{equation}\label{T}
T=\sum_{k=1}^d\left(Z^k-\frac{SY^k}{2\sqrt{m}}\right)\hat\e_k,
\quad
\mbox{with} \quad SY^k=\frac{1}{\sqrt{m}}\sum_{i=1}^mY_i^k,
 \quad
\hat\e_k=\1_{\{SY^k> \sqrt{2\log d}\}},
\end{equation}
and consider the classifier $\psi=\1_{\{T>0\}}$. We can write
$SY^k=\e_k \lambda +\zeta_k$ where $\zeta_k$ are independent
standard normal random variables. It is well known that
$\max_{k=1,\dots,d}|\zeta_k|\le \sqrt{2\log d}$ with probability
tending to 1 as $d\to+\infty$. This and (\ref{1}) imply that, with
probability tending to 1, the vector $(\hat\e_1,\dots,\hat\e_d)$
recovers exactly $(\e_1,\dots,\e_d)$ and the statistic $T$ coincides
with
$$
\hat T=\sum_{k=1 : \e_k=1}^d
\left(Z^k-\frac{SY^k}{2\sqrt{m}}\right).
$$
Since $E^{(u)}(SY^k/\sqrt{m})=\e_k a_d$ and
${\rm{Var}}^{(u)}(SY^k/\sqrt{m})=1/m$, we find:
\begin{eqnarray*}
& & E^{(u)}_{H_0}(\hat T)=-\frac{a_d}{2}\sum_{k=1}^d\e_k,\quad
E^{(u)}_{H_1}(\hat T)=\frac{a_d}{2}\sum_{k=1}^d\e_k, \\
& & {\rm{Var}}^{(u)}_{H_0}(\hat T)={\rm{Var}}^{(u)}_{H_1}(\hat
T)=\left(1+\frac{1}{4m}\right)\sum_{k=1}^d\e_k.
\end{eqnarray*}
It follows from Chebyshev's inequality that under (\ref{0}) we have
\begin{equation}\label{lab}
E^{(u)}_{H_0}(\psi)=P^{(u)}_{H_0}(T>0)\to 0,\quad
E^{(u)}_{H_1}(1-\psi)=P^{(u)}_{H_1}(T\le 0)\to 0
\end{equation}
as $d\to +\infty$. Note that this argument is applicable in the
general model of Section 1 (since the convergence in (\ref{lab}) is
uniform in $u\in U_{\beta,a_d}$), implying successful classification
by $\psi$ under conditions (\ref{0}) and (\ref{1}).

3$^\circ$. Let us now discuss a connection between conditions
(\ref{0}) and (\ref{1}). First, (\ref{1}) implies (\ref{0}) if $m$
is not too large:
\begin{equation}\label{2a}
m=o\left(d^{1-\beta}\log d\right).
\end{equation}
On the other hand, if $m$ is very large:
\begin{equation}\label{2b}
\exists \ b>0:\ \ m\ge b d^{1-\beta}\log d,
\end{equation}
then we have $x^2\ge ba_d^2d^{1-\beta}$, and condition (\ref{0})
implies (\ref{1}). Thus, the relation
$$
d^{1-\beta}a_d^2\asymp 1
$$
determines the classification boundary in the general model of
Section 1 if $m$ is very large (satisfies (\ref{2b})).

4$^\circ$. Finally, note that we can control conditions (\ref{1})
and (\ref{1*}) by their data-driven counterparts. In fact,
$\max_{k=1,\dots,d}|\zeta_k|\le \sqrt{2\log d}$ with probability
tending to 1 as $d\to+\infty$. Hence, if (\ref{1*}) holds, then
$M_Y\eqdef \max_{1\le k\le d}SY^k\le 3\sqrt{2\log d}$ with the same
probability. It is therefore convenient to consider the following
pre-classifier taking values in $\{0,1,ND\}$ ($ND$ means ``No
Decision", i.e., we need to switch to some other classifier):
$$
\psi^{pre}=\begin{cases} 0 &{\rm if}\quad T\le 0,\
M_Y>3\sqrt{2\log d },\\
1  &{\rm if}\quad T>0,\ M_Y>3\sqrt{2\log d},\\
ND &{\rm if}\quad M_Y\le 3\sqrt{2\log d },
\end{cases}
$$
where $T$ is given by (\ref{T}).
%Suppose that we choose a classifier $\bar \psi$ if
%$\psi^{pre}$ takes the value $ND$. Then the combined classifier is
%$\tilde\psi = \psi^{pre}\1_{\psi^{pre}\ne ND} + \bar \psi
%\1_{\psi^{pre}= ND}$ and we have
%\begin{eqnarray*}
%&&E^{(u)}_{H_0}(\tilde\psi) \le
%E^{(u)}_{H_0}\l(\bar\psi\1_{\{\psi^{pre}=ND\}}\r) +
%P^{(u)}_{H_0}(\psi^{pre}=1)
%=
%E^{(u)}_{H_0}\l(\bar\psi\1_{\{\psi^{pre}=ND\}}\r)+o(1),\\
%&&E^{(u)}_{H_1}(1-\tilde\psi) \le
%E^{(u)}_{H_1}\l((1-\bar\psi)\1_{\{\psi^{pre}=ND\}}\r) + o(1).
%\end{eqnarray*}
%In other words, $\tilde\psi$ has asymptotically the same the risk
%behaviour as $\bar \psi$ considered
The argument in 2$^\circ$ implies that $\psi^{pre}$ classifies
successfully if $ND$ is not chosen. Under condition (\ref{1*}) the
pre-classifier chooses $ND$ with probability tending to 1
 and then we apply
one of the classifiers suggested below in this paper. We prove their
optimality under assumption (\ref{1*}).

The above remarks can be easily extended to the case of Gaussian
errors with known variance $\sigma^2>0$  by using the normalization
$Z^k/\sigma, SY^k/\sigma$.  Moreover, they extend to the case
%s of unknown variance and
of non-Gaussian errors under the Cram\'er condition and the
additional assumption $m/\log d\to+\infty$ (cf. Section \ref{Ext}).

%%%%%%%%%%%%%%%%%%%%%%%%%%%%%

\section{Classification boundary for moderately sparse
vectors}\label{sec2}

In this section we consider the case of moderately sparse vectors.
To simplify the notation, we set without loss of generality
$\sigma=1$. Assume that $R_d=d^{1/2-\b}a_d$ satisfies:
\begin{equation}\label{upper1}
\lim_{d\to +\infty}R_d= +\infty
\end{equation}
and consider the  classifier based on a linear statistic:
$$
\psi^{lin}=\1_{\{T{'}>0\}},\quad T{'}=\sum_{k=1}^d
\left(Z^k-\frac{1}{2m}\sum_{i=1}^m Y_i^k\right).
$$
Note that $T{'}$ is similar to the statistic $T$ defined in
(\ref{T}) with the difference that in $T{'}$ we do not threshold to
estimate the positions of non-zero $\e_k$. Indeed, here we do not
necessarily assume (\ref{1}), and thus there is no guarantee that
$\e_k$ can be correctly recovered.

Assume that $\eta^k_j$ and $\xi^k_i$ for all $k,j,i$ are random
variables with zero mean and variance~1 (we do not suppose here that
$\eta^k_j$ have the same distribution as $\xi^k_i$). Then the means
of $Y_i^k$ and $Z^k$ are $E^{(u)}(Y_i^k)=E^{(u)}_{H_1}(Z^k)=\e_k
a_d$, $E^{(u)}_{H_0}(Z^k)=0$, their variances are equal to 1, and we
have:
\begin{eqnarray*}
& & E^{(u)}_{H_0}(T{'})=-\frac{a_d}{2}\sum_{k=1}^d\e_k,\quad
E^{(u)}_{H_1}(T{'})=\frac{a_d}{2}\sum_{k=1}^d\e_k, \\
& &  {\rm{Var}}^{(u)}_{H_0}(T{'})=
{\rm{Var}}^{(u)}_{H_1}(T{'})=d\left(1+\frac{1}{4m}\right).
\end{eqnarray*}
We consider now a vector $u\in \R^d$ of the form
\begin{equation}\label{upx}
u=(u_1,\dots,u_d):\ u_k=a_d\e_k ,\quad \e_k\in \{0,1\},\quad
\sum_{k=1}^d\e_k \ge cd^{1-\b}.
\end{equation}
By (\ref{upx}), Chebyshev's inequality and (\ref{upper1}), we obtain
\begin{eqnarray*}
E^{(u)}_{H_0}(\psi)= P^{(u)}_{H_0}(T'>0) & \le & P^{(u)}_{H_0}\Big(T'-
E^{(u)}_{H_0}(T{'})>cd^{1-\b}a_d/2\Big) \\
&  \le &
\frac{4d}{(cd^{1-\b}a_d)^2}\left(1+\frac{1}{4m}\right) \to 0
\end{eqnarray*}
as $d\to +\infty$. An analogous argument yields that $
E^{(u)}_{H_1}(1-\psi)\to 0$. The convergence here is uniform in $u$
satisfying (\ref{upx}), and thus uniform in $u\in U_{\beta, a_d}$.
Therefore, we have the following result.

%%%%%%%%
\begin{theorem}\label{uptr}
Let $\eta^k_j$ and $\xi^k_i$ for all $k,j,i$ be random variables
with zero mean and variance~1. If \nref{upper1} holds, then
successful classification is possible and it is realized by the
classifier $\psi^{lin}$.
\end{theorem}

%\begin{corollary}
%Let \begin{equation}\label{Rdup} m\asymp d^\g,\quad
%a_d=x/\sqrt{m},\quad \b\in (0,(1-\g)/2],\quad
%xd^{-\b+(1-\g)/2}\to\infty.
%\end{equation}
%Then successful classification is possible and it is realized by the
%classifier $\psi^{lin}$.
%\end{corollary}

\begin{remark}\label{R0}
{\rm We have proved theorem \ref{uptr} with the set of vectors $u$
defined by (\ref{upx}), which is larger than $U_{\beta, a_d}$. The
upper bound on $\sum_k \e_k$ in the definition of $U_{\beta, a_d}$
is not needed. Also the $\eta^k_j$ need not have the same
distribution as the $\xi^k_i$ and their variances need not be equal
to 1. It is easy to see that the result of theorem \ref{uptr}
remains valid if these random variables have unknown variances
uniformly bounded by an (unknown) constant.}
\end{remark}
%\begin{remark}\label{R0} {\rm It is straightforward to extend
%theorem \ref{uptr} to the noise $\eta^k_i$ with general variance
%$\sigma^2>0$. It suffices to replace $Z^k, Y_i^k$ by $Z^k/\sigma,
%Y_i^k/\sigma$ in the definition of the statistic, and to replace
%$a_d$ by $a_d/\sigma$ in the definition of $R_{d}$.
%An extension to unknown
%$\sigma^2$ is also possible (cf. Section \ref{ext1}).
% }
%\end{remark}

 The corresponding lower bound is given in the next theorem. For
 $a>0$, $t\in\R$, set
 $$
\ell_a(t)=f(t-a)/f(t),\quad D_a=\int\ell^2_a(t)f(t)dt,
$$
and
$$
 D_d(m,a,\b)=d^{1-2\b}D_a^m(D_a-1).
$$
\begin{theorem}\label{Lo1}
Let either $m\ge 1$ be fixed or $m=m_d\to +\infty$. If
\begin{equation}\label{lower.a}
\lim_{d\to +\infty}D_d(m,a_d,\b)= 0,
\end{equation}
 then successful classification
is impossible.
\end{theorem}
{\bf Proof of theorem \ref{Lo1}} is given in Section \ref{L}.

\begin{corollary}
Let $f$ be the density of standard normal distribution. If
\begin{equation}\label{Rd}
\lim_{d\to +\infty}R_d= 0,
\end{equation}
then successful classification is impossible for $\b\in (0,1/2]$ and
$m$ fixed or for $\b\in (0,1/2)$ and $m=m_d\to+\infty$ such that
$m=O(d^{1-2\b})$.
\end{corollary}
{\bf Proof}. For the standard normal errors we have $D_a=e^{a^2}$.
Therefore, condition \nref{lower.a} can be satisfied only if
$ma_d^2=o(1)$ as $d\to +\infty$. Moreover, in this case
\begin{equation}\label{DR}
 D_d(m,a_d,\b)\asymp
d^{1-2\b}a_d^2(1+ma_d^2)\asymp R_d^2.
\end{equation}
Thus, if $ma_d^2=o(1)$, conditions \nref{lower.a} and \nref{Rd} are
equivalent. Now, \nref{Rd} and the assumption $\b\in (0,1/2]$ imply
$a_d=o(1)$. This proves the corollary for fixed~$m$. Also, if $\b\in
(0,1/2)$ and $m=m_d\to+\infty$ such that $m=O(d^{1-2\b})$, then
$ma_d^2=O(R_d^2)=o(1)$.
\endproof

%\begin{corollary}
%Let $f$ be the density of standard normal distribution and
%\begin{equation}\label{Rdl}
%m\asymp d^\g,\quad a_d=x/\sqrt{m},\quad \b\in (0,(1-\g)/2],\quad
%x=o(d^{\b-(1-\g)/2}).
%\end{equation}
%Then successful classification is impossible.
%\end{corollary}
%{\bf Proof} follows directly to the proof above.

\begin{remark}\label{R2} {\rm
Relation \nref{DR} is valid for a larger class of noise
distributions, e.g., for non-Gaussian noise with finite Fisher
information. Indeed, assume that $\ell_a(t)$ is
$L_2(f)$-differentiable at point $a=0$, i.e., there exists a
function $\ell^\prime(\cdot)$ such that
\begin{equation}\label{diff}
\|\ell_a(\cdot)-1-a\ell^\prime(\cdot) \|_f=o(a),\quad
0<\|\ell^\prime(\cdot)\|_f<+\infty,
\end{equation}
where $\| g(\cdot) \|_f^2 = \int_{\R} g^2(x) \; f(x) \; dx$. Observe
that
$$
\|\ell^\prime(\cdot) \|_f^2=\int_\R\frac{(f^\prime(x))^2}{f(x)}dx
\eqdef I(f)
$$
is the Fisher information of $f$ (with $f'$ defined in a somewhat
stronger sense than, for instance, in \cite{IH}.
 Under assumption
(\ref{diff}) we have
\begin{equation*}\label{Lext2}
D_a=1+\|\ell_a(\cdot)-1\|_f^2,\
\|\ell_a(\cdot)-1\|_f^2=a^2(I(f)+o(1))
\end{equation*}
as $a\to 0$.
}
\end{remark}

Combining remarks \ref{R0} and \ref{R2} with theorems \ref{uptr} and
\ref{Lo1} we see that relation \nref{bound1} determines the
classification boundary for $\b\in (0,1/2]$ and fixed $m$ or for
$\b\in (0,1/2)$ and $m\to+\infty,\ m=O(d^{1-2\b})$, if the errors
have zero mean, finite variance and finite Fisher information.

\medskip

As corollaries of theorems \ref{uptr} and \ref{Lo1} we can establish
classification boundaries for particular choices of $a_d$. Recall
that non-trivial results can be expected only if $a_d$ satisfies
(\ref{1*}). For instance, consider $a_d=d^{-s}$ with some $s>0$.
Then for fixed $m$ the classification boundary in the region
$\b\in(0,1/2]$ is given by $s=\b-1/2$, i.e., successful
classification is possible if $s<1/2-\b$, and is impossible if
$s>1/2-\b$. Other choices of $a_d$ appear to be less interesting
when $\b\in(0,1/2]$. For example, in the next section we consider
the sequence $a_d = s\sigma\sqrt{(\log d)/m}$ with some $s>0$. If
$a_d$ is chosen in this way, successful classification is possible
for all $\b\in(0,1/2]$ with no exception, so that there is no
classification boundary in this range of $\b$.

Finally, note that theorem \ref{uptr} is valid for all $\beta\in
(0,1)$. However, for $\b> 1/2$ its assumption $\lim_{d\to
+\infty}R_d= +\infty$ guaranteeing successful classification is much
too restrictive as compared to the correct classification boundary
that we shall derive in the next section. The lower bound of theorem
\ref{Lo1} is also valid for all $\beta\in (0,1)$. However, we shall
see in the next section that it is not tight for highly sparse
vectors when $\b>3/4$ (cf. proof of theorem \ref{Lo2}).

%%%%%%%%%%%%%%%%%%%%%%%%%

\section{Classification boundary for highly sparse vectors}
We now analyse the case of highly sparse vectors, i.e.,we suppose
that $\b\in (1/2,1)$ if $\log m=o(\log d)$, and $\b^*=\b/(1-\g)\in
(1/2,1)$ if $\log m\sim \gamma \log d,\ \g\in (0,1)$. We shall show
that the classification boundary for this case is expressed in terms
of the function
\begin{equation*}\label{phii}
\phi(\b)= \begin{cases}
\phi_1(\b) & \rm{if} \ 1/2 < \b\le 3/4,\\
\phi_2(\b) & \rm{if} \ 3/4 < \b < 1,
\end{cases}
\end{equation*}
where the functions $\phi_1$ and $\phi_2$ are defined in
(\ref{phi12}). Note that $\phi_1$ and $\phi_2$ are monotone
increasing on $(1/2,1)$, satisfy $\phi_1(\b)\le \phi_2(\b)$ for all
$\b\in (1/2,1),$ and the equality $\phi_1(\b)=
\phi_2(\b)(=1/\sqrt{2})$ holds if and only if $\b=3/4$.

 The following notation will be useful in the sequel:
\begin{equation}\label{**}
T_d = \sqrt{\log d},\quad s=s_d = a_d/\sigma T_d,
\end{equation}
and
\begin{equation}\label{x}
x=s\sqrt{m},\quad x_0=sm/\sqrt{m+1},\quad x_1=s\sqrt{m+1},\quad
x^*=\frac{x}{1-\g}.
\end{equation}
Clearly, $x_0< x<x_1$. We allow $s, x, x_0, x_1$ to depend on $d$
but do not indicate this dependence in the notation for the sake of
brevity. We shall also suppose throughout that \nref{1*} holds, so
that $x_1=O(1)$ as $d\to+\infty$.

\subsection{Lower bound}\label{sec3l}
The next theorem gives a lower bound of classification for highly
sparse vectors.
\begin{theorem}\label{Lo2}
Let the noise density $f$ be Gaussian $\CN(0,\sigma^2)$,
$\sigma^2>0$. Assume that $\b\in (1/2,1)$ and
%$x_1<\phi(\b)$
$\limsup_{d\to+\infty}x_1<\phi(\b)$. Then successful classification
is impossible  for fixed $m$ and for $m=m_d\to +\infty$.
\end{theorem}
{\bf Proof of theorem \ref{Lo2}} is given in Section \ref{L}.

\smallskip

Though theorem \ref{Lo2} is valid with no restriction on $m$, it
does not provide a correct classification boundary if $m$ is large,
i.e., $\log m\sim \gamma \log d,\ \g\in (0,1)$, as in Scenarios (C)
and (D). The correct lower bound for large $m$ is given in the next
theorem.

\begin{theorem}\label{Lo2large}
Consider Scenario (C) with $\b^*=\b/(1-\g)\in (1/2,1)$ and
$$
a_d=\sigma x\sqrt{(\log d)/m}.
$$
Assume that $\limsup_{d\to+\infty}x^*<\phi(\b^*)$. Then successful
classification is impossible.
\end{theorem}
{\bf Proof of theorem \ref{Lo2large}} is given in Section \ref{L}.

\smallskip

Recall that, by an elementary argument, under Scenario (C) and for
$a_d$ as in theorem \ref{Lo2large}, successful classification is
impossible if $\b>1-\g$ (cf. remark after \nref{0}).  This is the
reason why in theorem \ref{Lo2large} we consider only $\b<1-\gamma$.

\subsection{Upper bounds for fixed $m$}\label{sec3}

We now propose optimal classifiers attaining the lower bound of
theorem \ref{Lo2} under Scenario~(A). First, we consider a procedure
that attains the classification boundary only for $\b\in [3/4,1)$
but has a simple structure. Introduce the statistics
$$
M_0=\max_{1\le k\le d}SY^k,\quad M=\max_{1\le k\le d}SZ^k
$$
where
\begin{equation}\label{SY}
SY^k=\frac{1}{\sqrt{m}}\sum_{i=1}^mY_i^k,\quad
SZ^k=\frac{1}{\sqrt{m+1}}\l(Z^k+\sum_{i=1}^mY_i^k\r).
\end{equation}
Define
$$
\Lambda_M=\frac{M}{\max(\sqrt{2}\,\sigma T_d, M_0)}.
$$
Taking a small $c_0>0$, consider the classifier of the form:
$$
\psi^{max}=\1_{\{\Lambda_M>1+c_0\}}.
$$
\begin{theorem}\label{Up1} Consider Scenario (A).
Let $\b\in (0,1)$ and \nref{1*} hold. Then, for any $c_0>0$,
\begin{equation}\label{up11}
\lim_{d\to +\infty} \sup_{u\in U_{\b,a_d}} E_{H_0}^{(u)} (\psi^{max})
= 0.
\end{equation}
If $\limsup_{d\to+\infty}x_1<\phi_2(\b)$,
%$x_1<\phi_2(\b)$
then, for any $c_0>0$,
\begin{equation}\label{up12}
\lim_{d\to +\infty} \sup_{u\in U_{\b,a_d}} E_{H_1}^{(u)} (\psi^{max})
= 0.
\end{equation}
If $\liminf_{d\to+\infty}x_1>\phi_2(\b)$
%$x_1>\phi_2(\b)$
, then there exists
$c_0>0$ such that
\begin{equation}\label{up13}
\lim_{d\to +\infty} \sup_{u\in U_{\b,a_d}} E_{H_1}^{(u)}
(1-\psi^{max}) = 0.
\end{equation}
\end{theorem}
{\bf Proof of theorem \ref{Up1}} is given in Section \ref{sec6}.

\bigskip

Theorems \ref{Lo2} and \ref{Up1} (cf. (\ref{up11}) and (\ref{up13})
and the fact that $\phi(\b)=\phi_2(\b)$ for $\b\in [3/4,1)$ ) imply
that $\psi^{max}$ attains the classification boundary for $\b\in
[3/4,1)$. On the other hand, (\ref{up12}) implies that for $\b\in
(1/2,3/4)$ (where $\phi(\b)=\phi_1(\b)<\phi_2(\b)$) the classifier
$\psi^{max}$ does not do the correct job. Its maximal risk ${\cal
R}_{\cal M}$ is asymptotically 1, which is larger than the risk 1/2
of the simple random guess. We therefore introduce another
classifier that has, however, a more involved structure. Consider
the statistics
\begin{eqnarray*}
L_0(t)&=&\sum_{k=1}^d(\1_{\{SY^k>t\sigma T_d\}}-\Phi(-tT_d)),\quad
\Delta_0(t)=\frac{L_0(t)}{\sqrt{d\Phi(-tT_d)}},\\
L(t)&=&\sum_{k=1}^d(\1_{\{SZ^k>t\sigma T_d\}}-\Phi(-t T_d)),\quad
\Delta(t)=\frac{L(t)}{\sqrt{d\Phi(-tT_d)}}
\end{eqnarray*}
where $t\in \R$, $\Phi$ is the standard normal cumulative
distribution function and the statistics $SY^k,\ SZ^k$ are defined
in \nref{SY}. Consider the grid
\begin{equation}\label{grid}
t_l=l h,\quad l=1,...,N,\quad t_N=\sqrt{2}\,\sigma,
\end{equation}
with a step $h>0$ depending on $d$ and such that $h=o(1)$, $T_d h\to
+\infty$. This implies that $1\ll N\ll T_d$ as $d\to +\infty$ (here
and below $v_d\ll w_d$ for $v_d>0$ and $w_d>0$ depending on $d$
means that $\lim_{d\to +\infty}v_d/w_d=0$). Set
$$
\Delta_0=\max_{1\le l\le N}\Delta_0(t_l),\quad \Delta=\max_{1\le
l\le N}\Delta(t_l),\quad \Lambda^*=\frac{\Delta}{H+\Delta_0},
$$
where $H=H_d$ is such that
\begin{equation}\label{H}
d^{bh}\ll H\ll d^B
\end{equation}
for any $B>0,\ b>0$ and any $d>d_0(B,b)$ where $d_0(B,b)$ is a
constant depending only on $B$ and $b$ (such an $H$ can be always
determined depending on the choice of $h$). Consider now the
classifier of the form
$$
\psi^*_m=\1_{\{\Lambda^*>H\}}.
$$
%Set $x_0=ms/\sqrt{m+1}=x_1m/(m+1)$.

\begin{theorem}\label{Up2}
Consider Scenario (A) with $\b\in (1/2,1)$ and assume \nref{1*}.
Then
\begin{equation}\label{up21}
\lim_{d\to +\infty} \sup_{u\in U_{\b,a_d}} E_{H_0}^{(u)} (\psi^*_m) =
0.
\end{equation}
If $\liminf_{d\to+\infty}x_1>\phi(\b)$ and
$\limsup_{d\to+\infty}x_0<\sqrt{2}$, then
\begin{equation}\label{up23}
\lim_{d\to +\infty} \sup_{u\in U_{\b,a_d}} E_{H_1}^{(u)} (1-\psi^*_m)
= 0.
\end{equation}
\end{theorem}
{\bf Proof of theorem \ref{Up2}} is given in Section \ref{sec6}.

\bigskip

Theorems \ref{Lo2}, \ref{Up1} and \ref{Up2} show that the
classification boundary for highly sparse vectors (i.e., for $\b\in
(1/2,1)$) is given by (\ref{detb}). Furthermore, the classifier
$\psi^*_m$ is optimal (attains the classification boundary) for
$\b\in (1/2,1)$, except for the case
$\limsup_{d\to+\infty}x_0\ge\sqrt{2}$, which is already covered by
the classifier $\psi^{max}$. Indeed, $x_0\ge\sqrt{2}$ implies that
$x_1\ge \sqrt{2}(1+1/m)>\phi_2(\beta)$ for all $\b\in (1/2,1)$.

%%%%%%%%%%%%%%%

\subsection{Upper bounds for
$m\to+\infty,\ \log m=o(\log d)$} \label{sec4}

In this subsection we analyse Scenario (B). Then $m=m_d\to +\infty$,
$\log m=o(\log d)$ as $d\to +\infty$ and the classifier $\psi^*_m$
is not, in general, optimal. Nevertheless, we propose another
classifier $\psi^*_\infty$, which attains essentially the same
classification boundary as in Subsection \ref{sec3} above. Introduce
the statistics
$$
\Delta(t)=\frac{1}{\sigma\sqrt{d \Phi(-tT_d)}}
\sum_{k=1}^dZ^k\1_{\{SY^k>t\sigma T_d\}},\quad \Delta=\max_{1\le
l\le N}\Delta(t_l),
$$
where the maximum is taken over the grid \nref{grid}. Here and below
we use the same notation $\Delta(t)$, $\Delta$ as previously for
different ratio statistics, since it causes no ambiguity. Set also
$$
\Delta_*=\sum_{k=1}^d\1_{\{SY^k>\sqrt{2}\,\sigma T_d\}}
$$
and define
\begin{equation}\label{upperM.3}
\Lambda^*_\infty=\frac{\Delta}{\sqrt{H+\Delta_*}},\quad
\psi^*_\infty=\1_{\{\Lambda^*_\infty>H\}},
\end{equation}
where $H$ satisfies \nref{H}.

\begin{theorem}\label{Up3}
Consider Scenario (B). Let $\b\in (1/2,1)$  and let \nref{1*} hold.
Then
\begin{equation}\label{up31}
\lim_{d\to +\infty} \sup_{u\in U_{\b,a_d}} E_{H_0}^{(u)}
(\psi^*_\infty) = 0.
\end{equation}
If %$a_d \sqrt{m_d/\log d}>\phi(\b)$ for all $d$,
$\liminf_{d\to+\infty}x>\phi(\b)$, then
\begin{equation}\label{up33}
\lim_{d\to +\infty} \sup_{u\in U_{\b,a_d}} E_{H_1}^{(u)}
(1-\psi^*_\infty) = 0.
\end{equation}
\end{theorem}
{\bf Proof of theorem \ref{Up3}} is given in Section \ref{sec6}.

%%%%%%%%%%%%%%%%%%%%%%%%%%%

\subsection{Upper bound for Scenario (C)} \label{sec5}

We now suggest an asymptotically optimal classifier for Scenario
(C). For $t\ge 0$ we introduce the statistics
$$
L^1(t)=\sum_{k=1}^d Z^k\1_{\{SY^k>\sigma tT_d\}},\quad
L^0(t)=\sum_{k=1}^d \1_{\{SY^k>\sigma tT_d\}},\quad
\Delta(t)=\frac{L^1(t)}{\sigma \sqrt{N^2+L^0(t)}},
$$
Take a grid $t_1,\dots,t_N$ of the form \nref{grid} and define the
classifier
$$
\psi_\infty=\1_{\{\Delta>\,4\,N\}},\quad \mbox{where}\quad
\Delta=\max_{1\le l\le N}\Delta(t_l).
$$

\begin{theorem}\label{Up4} Consider Scenario (C) with
$ a_d=\sigma x\sqrt{(\log d)/m}$. Let $\b^*=\b/(1-\g)\in (1/2,1)$
and let \nref{1*} hold. Then
\begin{equation}\label{up31large}
\lim_{d\to +\infty} \sup_{u\in U_{\b,a_d}} E_{H_0}^{(u)}
(\psi_\infty) = 0.
\end{equation}
If $\liminf_{d\to+\infty}x^*>\phi(\b^*)$, then
\begin{equation}\label{up33large}
\lim_{d\to +\infty} \sup_{u\in U_{\b,a_d}} E_{H_1}^{(u)}
(1-\psi_\infty) = 0.
\end{equation}
\end{theorem}
{\bf Proof of theorem \ref{Up4}} is given in Section \ref{sec6}.
%\begin{remark}\label{R.up}{\rm
%The analyze of the proof of theorem \ref{Up4} shows the the
%statements of Theorem hold true as $m\to\infty,\ \log(m)\sim\g\log
%(d)$ with $0\le \g<1$. In particular the classifier $\psi_\infty$
%provides successful classification under assumptions of Theorem
%\ref{Up3} as well.}
%\end{remark}
%%%%%%%%%%%%%%%

\section{Extensions}\label{Ext}

\subsection{Unknown variances}\label{ext1}

The classifiers proposed in the previous section can be easily
extended to the model with unknown variance $\sigma^2$, so that the
results of theorems \ref{Up1}, \ref{Up2}, \ref{Up3} and \ref{Up4}
remain valid. We present here the general lines of such a
modification without going into the details of the proofs that do
not differ much from those in Section \ref{sec6}.

First, note that there exists an estimator $\hat\sigma^2_d$
satisfying
\begin{eqnarray}\label{sigma}
\hat\sigma_d^2&=&\sigma^2+\eta_d,\
\end{eqnarray}
where $\eta_d\to 0$ in $P_{H_0}^{(u)}$-probability, and
\begin{eqnarray}\label{sigma1}
\hat\sigma_d^2&=&\sigma^2+O(d^{-\b}a_d^2)+
(1+d^{-\b/2}a_d)^{1/2}\eta_d',\
\end{eqnarray}
where $\eta_d'\to 0$ in $P_{H_1}^{(u)}$-probability, uniformly in
$u\in U_{\b,a_d}$, as $d\to +\infty$.

For example, we can take the standard sample variance
$$ \hat\sigma^2_d=\frac
1d \sum_{k=1}^d(Z^k)^2.
$$
Assume that $\eta_j^k$ are i.i.d. $\CN(0,\sigma^2)$ random variables
with unknown $\sigma$. Then (\ref{sigma}) and (\ref{sigma1}) are
satisfied. In fact, it is easy to see that
$$
E_{H_0}^{(u)}(\hat\sigma_d^2)= \sigma^2,\quad
E_{H_1}^{(u)}(\hat\sigma_d^2)= \sigma^2 +\frac 1d \sum_{k=1}^d
u_k^2= \sigma^2+O(a_d^2d^{-\b}),
$$
and analogously
$$
\Var_{H_0}^{(u)}(\hat\sigma_d^2)=\frac {2\sigma^4}{d},\quad
\Var_{H_1}^{(u)}(\hat\sigma_d^2) =
\frac{1}{d}\l(2\sigma^4+O(d^{-\b}a_d^2)\r)=o(1+d^{-\b}a_d^2)
$$
as $d\to +\infty$. Applying Chebyshev's inequality, we get
\nref{sigma} and \nref{sigma1}. We also note that these relations
hold under much weaker assumptions than the normality of $\eta_j^k$.
It suffices to have, for example, independent random variables
$\eta_j^k$ such that $E(\eta_j^k)=0$, $E[(\eta_j^k)^2]=\sigma^2$ and
$\max_{j,k} E[(\eta_j^k)^4]<+\infty$.

We now discuss how to modify the proposed classifiers using
$\hat\sigma_d$. For $\psi^{pre}$ and $\psi^{max}$, we replace the
unknown $\sigma$ in their definitions by $\hat\sigma_d$ and change
$\sqrt{2\log d}$ into $\sqrt{b\log d},\ b>2$ for $\psi^{pre}$. If
$R_{d}=O(1)$ (which is the case for highly sparse vectors under
(\ref{**})), then $d^{-\b}a_d^2=o(1)$ and (\ref{sigma}) implies that
the ratio $\hat\sigma_d/\sigma$ is close to 1 in
$P_{H_1}^{(u)}$-probability as well. Therefore, for the study of the
variance modified versions of classifiers $\psi^{pre}, \psi^{max}$,
we can use not only (\ref{sigma}) but also the fact that
$\hat\sigma_d^2=\sigma^2+\bar\eta_d$ where $\bar\eta_d\to 0$ in
$P_{H_1}^{(u)}$-probability. Thus, the desired upper bounds for
these classifiers follow in an easy way from the results in Section
\ref{sec3}.

For the classifier $\psi^*_m$, we replace the statistics $L_0(t),
L(t), \Delta_0(t), \Delta(t)$ by
\begin{eqnarray*}
L_0(t)&=&\sum_{k=1}^d(\1_{\{SY^k>tT_d\}}-\1_{\{Z^k>tT_d\}}),\quad
\Delta_0(t)=\frac{L_0(t)}{\sqrt{d\Phi(- tT_d/\hat\sigma_d)}},\\
L(t)&=&\sum_{k=1}^d(\1_{\{SZ^k>tT_d\}}-\1_{\{Z^k>tT_d\}}),\quad
\Delta(t)=\frac{L(t)}{\sqrt{d\Phi(-tT_d/\hat\sigma_d)}},
\end{eqnarray*}
and we take a grid
$$
t_l=l h,\ l=1,...,N, t_N=\sqrt{2}\,\hat\sigma_d+O(h),
$$
with step $h$ as in \nref{grid}. The cardinality $N$ of the grid
thus becomes a random variable. However, the relation $N=O(T_d)$
holds true in probability under \nref{sigma}. Note that the modified
statistics $\Delta_0(t)$ and $\Delta(t)$ contain the additional
factor $A(t)=\sqrt{\Phi(-tT_d/\sigma)/\Phi(-tT_d/\hat\sigma_d)}$ as
compared to the original ones. If \nref{sigma} holds, these factors
are (in probability) of the form $\exp(o(T_d))$ uniformly in
$t=O(1)$. Under $P^{(u)}_{H_s},\, s=0,1$, the expectations of the
summands with $\e_k=0$ in $L_0(t)$ and $L(t)$ vanish. The other
elements of the proof for the modified statistics are similar to
those in Section \ref{sec6}.

For the classifier $\psi^*_\infty$, we replace the statistics $
\Delta(t)$ and $\Delta_0$ by
$$
\Delta(t)=\frac{1}{\hat\sigma_d\sqrt{d\Phi(-tT_d/\hat\sigma_d)}}
\sum_{k=1}^dZ^k\1_{\{SY^k>tT_d\}},\quad \Delta=\max_{1\le l\le
N}\Delta(t_l),
$$
with the same grid as above, and
$$
\Delta_0=\sum_{k=1}^d\1_{\{SY^k>\hat\sigma_d\sqrt{2}T_d\}}.
$$
We make similar modifications for the classifier $\psi_\infty$. The
arguments above are enough for the proof of Sections \ref{U2},
\ref{U4} to hold through.

\subsection{Non-Gaussian noise}\label{ext2}

We now discuss an extension of our results to Scenario (D). The
remarks on the pre-classifier $\psi^{pre}$ in Section \ref{pre} and
the proofs of the upper bounds in Section \ref{sec6} are only based
on the constraint \nref{1*} and the following property of the tails
of the Gaussian distribution:
\begin{equation}\label{GT}
\log P(S\zeta^m>\sigma t)\sim -\frac{t^2}{2}, \quad t \in
[U_0,U_1]\quad \text{for $U_0\to+\infty$ and $U_1=O(T_d)$.}
\end{equation}
Here $S\zeta^m=\frac{1}{\sqrt{m}}\sum_{i=1}^m\zeta_i$, and $\zeta_i$
are i.i.d.  ${\cal N}(0,\sigma^2)$ random variables.

Indeed, in Subsection \ref{U4}  we can write
$\Phi(-tT_d)$ as $P(S\zeta^m>\sigma tT_d)$. From (\ref{GT}) we
deduce
$$
P(S\zeta^m>tT_d/\sigma)=A_d d^{-(t_+)^2/2},\quad t_+=\max(0,t),
$$
where $A_d$ satisfies \nref{AA} for $t_+=O(1)$. This is exactly the
relation \nref{A}, which is also the only property of the noise
distribution needed for the proofs in Subsection
\ref{U4}.

If $m$ is large enough, relation (\ref{GT}) holds not only for the
Gaussian $\zeta_i$. It suffices to have the i.i.d. $\zeta_i$ with
$E\zeta_i=0,\ E(\zeta_i^2)=\sigma^2>0$ satisfying the Cram\'er
condition:
$$ %\begin{equation}\label{cramer}
\exists\quad h_0>0 :\quad E \left( e^{h\zeta_i} \right) <+\infty,\quad \forall \
h\in(-h_0, h_0).
$$ %\end{equation}
and $m\gg \log d$. In fact, using theorem 5.23 in  \cite{Pe} we get
that, under the Cram\'er condition and for $t=o(\sqrt{m})$,
$$
P(S\zeta^m>\sigma
t)=\Phi(-t)\exp\l(\frac{t^3}{\sqrt{m}}\lambda\l(\frac{t}{\sqrt{m}}\r)\r)\l\{1+O\l(\frac{t+1}{m}\r)\r\},
$$
 where $\lambda(t)$ is the Cram\'er series. Inserting here the
expression for the Cram\'er series and
the relation
$\log \Phi(-t)=-t^2/2-\log t+O(1)$ as $t\to+\infty$, we obtain
$$
\log P(S\zeta^m>\sigma
t)=-\frac{t^2}{2}\l(1+O\l(\frac{t}{\sqrt{m}}\r)+o(1)\r)\sim
-\frac{t^2}{2}
$$
as $t\to+\infty,\ t=o(\sqrt{m})$. These remarks  allow us to follow
the proof of theorem \ref{Up4} in Section~\ref{sec6} leading to the
next result.
%%%%%%%%%%%%%%%
\begin{theorem}\label{Tup3}
Consider Scenario (D) with $ a_d=\sigma x\sqrt{(\log d)/m}$. Let
$\b^*=b/(1-\g)\in (1/2,1)$ and let \nref{1*} hold. Then
\begin{equation}\label{up31a}
\lim_{d\to +\infty} \sup_{u\in U_{\b,a_d}} E_{H_0}^{(u)}
(\psi_\infty) = 0.
\end{equation}
If $\liminf_{d \to +\infty} x^*>\phi(\b^*)$, then
\begin{equation}\label{up33a}
\lim_{d\to +\infty} \sup_{u\in U_{\b,a_d}} E_{H_1}^{(u)}
(1-\psi_\infty) = 0.
\end{equation}
\end{theorem}

%%%%%%%%%%%%%%%%%

\subsection{Adaptive procedures}\label{sec_ad}

 We have proposed several classifiers, which attain the classification
boundary under various conditions on $m, a_d, \b$. In order to
obtain an adaptive procedure that attains this boundary
simultaneously for several domains of $m, a_d, \b$, it suffices to
combine the classifiers in the following way. We start with the
pre-classifier $\psi^{pre}$. If it outputs ``No Decision", then we
combine the classifiers $\psi^{lin}$, $\psi^{max}$ and $\psi^*_m$
using the Bonferroni device, i.e., our classifier will be
$\max(\psi^{lin},\psi^{max},\psi^*_m)$. This means that we allocate
$Z$ to the $P_Y$-population iff it is allocated to $P_Y$ by at least
one of the three classifiers. Analogously, if $m\to+\infty$, then we
classify by $\max(\psi^{lin},\psi^*_\infty)$ or by
$\max(\psi^{lin},\psi_\infty)$.

\section{Proof of the lower bounds}\label{L}

In this section we prove theorems \ref{Lo1}, \ref{Lo2} and
\ref{Lo2large}. Without loss of generality we consider only the case
${\cal R}_{\max}=1/2$ (cf. (\ref{bi})). Observe that if a
probability measure $\mu^d$ on $\R^d$ is such that
\begin{eqnarray}\label{mu1}
\mu^d(U_{\b,a})=1+o(1)
\end{eqnarray}
 as $d\to +\infty$, then
\begin{eqnarray} \sup_{u\in U_{\b,a_d}}{\cal R}_M(\psi)&\ge& \max
\Big(\int E_{H_0}^{(u)} (\psi)\mu^d(du),\nonumber
\int E_{H_1}^{(u)}(1-\psi)\mu^d(du)\Big) + o(1)\\
&\ge& \frac{1}{2} \Big(\int E_{H_0}^{(u)} (\psi)\mu^d(du)+ \int
E_{H_1}^{(u)}(1-\psi)\mu^d(du)\Big) + o(1)\nonumber\\
&=& \frac{1}{2}\int \left(\psi + (1-\psi)\frac{d{\sf P}_{H_1}}
{d{\sf P}_{H_0}}\right) d{\sf P}_{H_0} + o(1)\label{bay}
\end{eqnarray}
where ${\sf P}_{H_s}$, $s=0,1$, are the ``posterior" probability
measures defined by
$$
{\sf P}_{H_s}(A) = \int P_{H_s}^{(u)}(A) \mu^d(du)
$$
for any Borel set $A$ of $\left(\R^d\right)^m \times \R^d$.
 In view of (\ref{bay})
%and of the Fatou lemma
, if the likelihood ratio
$$L({\bf Y},Z)=\frac{d{\sf P}_{H_1}} {d{\sf P}_{H_0}}$$ satisfies
\begin{equation}\label{aim}
L({\bf Y},Z) \to 1 \quad \text{in ${\sf P}_{H_0}$-probability}
\end{equation}
as $d\to +\infty$, then the left-hand side of (\ref{bi}) is greater
than or equal to 1/2. This immediately entails (\ref{bi}) because
the risk of the simple random guess classifier equals 1/2.

 Since ${\sf E}_{H_0}(L({\bf Y},Z)-1)^2={\sf
E}_{H_0}L^2({\bf Y},Z) -1$, relation (\ref{aim}) holds if
\begin{equation}\label{aim1}
\limsup_{d\to +\infty}{\sf E}_{H_0}L^2({\bf Y},Z) = 1.
\end{equation}
Based on these remarks, the proofs of theorems \ref{Lo1} and
\ref{Lo2} will proceed by constructing a prior measure $\mu^d$
satisfying (\ref{mu1}) and proving (\ref{aim1}).

In this section we assume without loss of generality  that
$\sigma=1$ and that the constants $c,C$ in the definition of
$U_{\b,a_d}$ are such that $c<1<C$. %We also set for brevity $a=a_d$.

The prior measure that we choose here is of the form
$\mu^d(du)=\prod_{k=1}^d\mu(du_k)$ where
$\mu=(1-p)\delta_0+p\delta_{a_d}$, $p=d^{-\b}$, and $\delta_t$ is the
Dirac mass at point $t\in\R$. In other words, the prior measure
corresponds to $u_k=a_d \e_k $ with i.i.d Bernoulli random entries $
\e_k$ that take value 1 with probability $p=d^{-\b}$ and value 0
with probability $1-d^{-\b}$.

\begin{lemma}\label{L1} Let $ 0<c<1<C<+\infty$. Then the prior
measure $\mu^d$ defined above satisfies (\ref{mu1}).
\end{lemma}
{\bf{Proof of  Lemma \ref{L1}}}. Set $G(u)=\sum_{k=1}^d u_k$. We
have to check that
$$
\mu^d(G(u)>a_d Cd^{1-\b})\to 0,\quad \mu^d(G(u)<a_d cd^{1-\b})\to 0.
$$
Since
$$
E_{\mu^d}(G(u))=a_d dp=a_d d^{1-\b},\quad
\Var_{\mu^d}(G(u))=a_d^2 dp(1-p)\sim a_d^2 d^{1-\b}\quad \forall\
\b\in(0,1),
$$
it follows from Chebyshev's inequality that
\begin{eqnarray*}
\mu^d(G(u)>a_d Cd^{1-\b})
%&=&
%\mu^d(G(u)-E_{\mu^d}(G(u))>a_d (C-1)d^{1-\b})\\
%\end{eqnarray*}
%\begin{eqnarray*}
&\le& \frac{1}{ d^{1-\b}(C-1)^2}\to 0, \\
%\end{eqnarray*}
%\begin{eqnarray*}
\mu^d(G(u)<a_d cd^{1-\b})
%&=&
%\mu^d(E_{\mu^d}(G(u))-G(u)>a_d (1-c)d^{1-\b})\\
%\end{eqnarray*}
%\begin{eqnarray*}
&\le& \frac{1}{d^{1-\b}(1-c)^2}\to 0
\end{eqnarray*}
as $d\to +\infty$.
\endproof

It remains now to prove that (\ref{aim1}) holds under
the assumptions of theorems \ref{Lo1} and \ref{Lo2}.

We shall need some notation. For $a \in \R$ define the probability
densities
\begin{equation}\label{gden}
f_a(\bY^k)=\prod_{i=1}^m f(Y_i^k-a), \quad
f_a(\bY^k,Z^k)=f_a(\bY^k)f(Z^k-a).
\end{equation}
%and
%$$ f_{a_d}(\bY^k)=\prod_{j=1}^m
%f(Y_i^k-a_d), \quad
%f_{a_d}(\bY^k,Z^k)=f_{a_d}(\bY^k)f(Z^k-a_d).
%$$
Let $P_0=\prod_{k=1}^dP_{0,k}$ be the probability measure that
corresponds to the pure noise. Here the measure $P_{0,k}$ has the
density $f_0(\bY^k,Z^k)=f_0(\bY^k)f(Z^k)$ and $E_{0,k}\left( \cdot
\right)$ denotes the expectation under $P_{0,k}$.

Next, write ${\sf P}_{H_s}=\prod_{k=1}^d P_{H_s,k}, \ s=0,1,$ where
the probability measures $P_{H_s,k}\ s=0,1,$ have the densities
$$ f_{H_s,k}(\bY^k, Z^k)=
(1-p)f_0(\bY^k) f(Z^k)+pf_{a_d}(\bY^k)f(Z^k-s a_d).
$$
%$$ f_{H_0,k}(\bY^k, Z^k)=
%f_0(Z^k)\l\{(1-p)f_0(\bY^k)+pf_{a_d}(\bY^k)\r\}
%$$
%and
%$$
%f_{H_1,k}(\bY^k, Z^k)=(1-p)f_0(\bY^k,Z^k)+pf_{a_d}(\bY^k,Z^k).
%$$
We denote by $E_{H_s,k},s=0,1$, the corresponding expectations. The
measures ${\sf P}_{H_s}$ have the following densities :
$$
f_{H_s}(\bY, Z)=\prod_{k=1}^d f_{H_s,k}(\bY^k, Z^k).
$$
 The
likelihood ratio is of the form
$$
L(\bY, Z)=\frac{d{\sf P}_{H_1}}{d{\sf
P}_{H_0}}=\prod_{k=1}^dL_k(\bY^k,Z^k),
$$
where
$$
L_k(\bY^k,Z^k)= \frac{(1-p)+pL(\bY^k,Z^k)}{(1-p)+pL(\bY^k)}=
\frac{1+p(L(\bY^k,Z^k)-1)}{1+p(L(\bY^k)-1)},
$$
and we set
$$
L(\bY^k)=\prod_{j=1}^m\ell_{a_d}(Y^k_j),\quad
L(\bY^k,Z^k)=L(\bY^k)\ell_{a_d}(Z^k)
$$
where $\ell_{a_d}(t)=f(t-a_d)/f(t)$. It will be convenient to
write $L_k$ in the form
\begin{equation}\label{lower00}
L_k(\bY^k,Z^k)=1+\Delta_k,\quad
\Delta_k=\frac{pL(\bY^k)(\ell_{a_d}(Z^k)-1)}{1+p(L(\bY^k)-1)}.
\end{equation}

%%%%%%%%%%%%%%%%%%%%%%

\subsection{Proof of theorem \ref{Lo1}}
Recall that
$
\frac{dP_{H_0,k}}{dP_{0,k}}=1+p(L(\bY^k)-1).
$
Since $E_{H_0,k}(\Delta_k)=0$, we obtain
\begin{eqnarray}\nonumber
{\sf E}_{H_0}(L^2(\bY,Z))&=&\prod_{k=1}^d\l(1+E_{H_0,k}\Delta_k^2\r)
%\\
%\nonumber
\le \exp\l(\sum_{k=1}^d E_{H_0,k}\Delta_k^2\r)\\
\nonumber &=& \exp\l(\sum_{k=1}^d E_{0,k}\l(\Delta_k^2
\frac{dP_{H_0,k}}{dP_{0,k}}\r)\r)\\ \nonumber
%&=&
%\exp\l(\sum_{k=1}^d
%E_{0,k}\l(\frac{p^2L^2(\bY^k)(\ell_{a_d}(Z^k)-1)^2}{1+p(L(\bY^k)-1)}\r)\r)\\
%\nonumber
&\le&\exp\l(\frac{p^2}{1-p}\,\sum_{k=1}^d
E_{0,k}L^2(\bY^k)\,
E_{0,k}(\ell_{a_d}(Z^k)-1)^2\r)\\
&=& \exp\l(\frac{dp^2D_{a_d}^m(D_{a_d}-1)}{1-p} \r),\label{lower2a}
\end{eqnarray}
where
$$
D_{a_d}=\int_{\R}\ell_{a_d}^2(t)f(t)dt=\int_{\R} \frac{f^2(t-a_d)}{f(t)}dt.
$$
Since $p=d^{-\b}\to 0$, relation
\nref{aim1} holds if
\begin{equation}\label{lower2}
D_d(m,a_d,\b)=d^{1-2\b}D_{a_d}^m(D_{a_d}-1)\to 0.
\end{equation}
This completes the proof of theorem \ref{Lo1}.

%%%%%%%%%%%%%%%%%%%%%%%%%

\subsection{Proof of theorem \ref{Lo2}}\label{ProofLo2}

Assume w.l.o.g. that $\sigma=1$.  Then $f$ is the standard normal
density, and thus $D_a=e^{a^2}$. We shall assume that $x_1$ is
fixed; the general case can be treated in a similar way by passing
to subsequences $x_{1,d}\to x_1\ge 0$.  By \nref{ad}, the condition
\nref{lower2} takes the form
\begin{equation}\label{lower6}
d^{1-2\b}\exp((m+1)a_d^2)=d^{1-2\b+x_1^2}\to 0,
\end{equation}
In other terms, the proof of theorem \ref{Lo1} implies that
successful classification is impossible if
\begin{equation}\label{lower7} x_1^2-2\b+1<0.
\end{equation}
This bound applies for any $\b\in(1/2,1)$, and it yields the result
of theorem \ref{Lo2} for $\b\in(1/2,3/4]$. It remains to show that a
bound better than \nref{lower7} can be obtained for $\b\in (3/4,1)$,
namely
\begin{equation}\label{lower8}
x_1^2-2\Big(1-\sqrt{1-\b}\Big)^2< 0.
\end{equation}
In order to prove this, set
$$
SY_{k}=\sum_{j=1}^m Y_j^k,\quad SZ_k=SY_k+Z^k,\quad k=1,...,d,\quad
T_{l,d}=\sqrt{2l\log d},
$$
and introduce the events
$$
\CA_{SY,k}=\{SY_k<T_{m,d}\} ,\quad \CA_{SZ,k}=\{SZ_k<T_{m+1,d}\},
$$
$$
\CA_{SY}=\bigcap_{k=1}^d\CA_{SY,k},\quad
\CA_{SZ}=\bigcap_{k=1}^d\CA_{SZ,k}.
$$
Observe that since
$$
P_{0,k}(SY_k\ge T_{m,d})=P_{0,k}(SZ_k\ge
T_{m+1,d})=\Phi\l(-\sqrt{2\log d}\r)=o(d^{-1}),
$$
we have
$$
P_0(\CA_{SY})\to 1,\quad P_0(\CA_{SZ})\to 1.
$$
Moreover
\begin{eqnarray*}
P_{H_0,k}(SY_k\ge T_{m,d})&=&P_{H_1,k}(SY_k\ge
T_{m,d}) \\
& = & (1-p)P_{0,k}(SY_k\ge
T_{m,d})+pP_{0,k}(SY_k\ge T_{m,d}-ma_d),\\
pP_{0,k}(SY_k\ge T_{m,d}-ma_d) & = & d^{-\b}\Phi\l(a_d
\sqrt{m}-\sqrt{2\log d}\r) \\
& < & d^{-\b}\Phi\l(a_d \sqrt{m+1}-\sqrt{2\log d}\r)
\asymp \frac{d^{-g}}{\sqrt{\log d}} =o(d^{-1}),
\end{eqnarray*}
where $ g \eqdef \b+\l(\sqrt{2}- x_1 \r)^2/2 \ge 1$ in view of
\nref{lower8}. Analogously, we have for $s=0,1$ :
\begin{eqnarray*}
P_{H_s,k}(SZ_k\ge T_{m+1,d})&=&(1-p)P_{0,k}(SZ_k\ge
T_{m+1,d})\\
&&+pP_{0,k}(SZ_k\ge T_{m+1,d}- (m+s) a_d),\\
pP_{0,k}(SZ_k\ge T_{m+1,d}-m a_d)&\le & pP_{0,k}(SZ_k\ge
T_{m+1,d}-(m+1) a_d)\\
&=&d^{-\b}\Phi\l(a_d \sqrt{m+1}-\sqrt{2\log d}\r)\\
&\asymp& \frac{d^{-g}}{\sqrt{\log d}}=o(d^{-1}).
\end{eqnarray*}
Thus,
\begin{equation}\label{lower9}
{\sf P}_{H_s}(\CA_{SY})\to 1,\quad {\sf P}_{H_s}(\CA_{SZ})\to
1,\quad s=0,1,
\end{equation}
as $d\to +\infty$. Set
$
\hat
L_k(\bY^k,Z^k)=L_k(\bY^k,Z^k)\1_{\{\CA_{SY,k}\bigcap\CA_{SZ,k}\}},\quad
\hat\Delta_k=\Delta_k\1_{\{\CA_{SY,k}\bigcap\CA_{SZ,k}\}},
$
where $\Delta_k$ is defined by \nref{lower00}, and
$
\hat L(\bY,Z)=\prod_{k=1}^d \hat L_k(\bY,Z).
$
Using (\ref{lower9}) we get that the main term in (\ref{bay})
satisfies
\begin{eqnarray*}
\int \left(\psi + (1-\psi) L(\bY,Z) \right) d{\sf P}_{H_0}
& = & \int\left(\psi + (1-\psi)\hat L(\bY,Z)\right) d{\sf P}_{H_0} \\
& & + \int (1-\psi)\1_{\{\bar\CA_{SY}\bigcup\bar\CA_{SZ}\}} d{\sf
P}_{H_1} \\
& = & \int\left(\psi + (1-\psi)\hat L(\bY,Z)\right) d{\sf P}_{H_0}
+o(1)
\end{eqnarray*}
as $d\to +\infty$. Repeating the argument after (\ref{bay}) we see
that to prove the theorem it suffices to show that
\begin{equation}\label{lower10}
\hat L(\bY,Z)\to 1 \quad \text{in ${\sf P}_{H_0}$-probability}.
\end{equation}
Using \nref{lower9} we obtain that $ {\sf E}_{H_0}\hat L(\bY,Z)={\sf
P}_{H_1}(\CA_{SY}\cap\CA_{SZ})\to 1. $ Therefore, to show
(\ref{lower10})  it suffices to prove that (cf. (\ref{aim1})):
\begin{equation}\label{lower10a}
\limsup_{d\to +\infty}{\sf E}_{H_0}\hat L^2(\bY,Z) = 1.
\end{equation}
We now prove (\ref{lower10a}). First note that, as follows from the
displays preceding (\ref{lower9}),
$$
E_{H_0,k}(\hat\Delta_k)=P_{H_1,k}(\CA_{SY,k}\cap\CA_{SZ,k})-
P_{H_0,k}(\CA_{SY,k}\cap\CA_{SZ,k})=o(d^{-1}),
$$
and
$$
0\le \hat L_k(\bY^k,Z^k)\le 1+\hat\Delta_k.
$$
Therefore, arguing as in \nref{lower2a} we obtain
\begin{eqnarray}\nonumber
{\sf E}_{H_0}(\hat L^2(\bY,Z))&=&\prod_{k=1}^d  E_{H_0,k}(\hat
L^2_k(\bY^k,Z^k))\le
\prod_{k=1}^d\l(1+E_{H_0,k}\hat \Delta_k^2+2E_{H_0,k}\hat \Delta_k\r)\\
\nonumber &\le& \exp\l(\sum_{k=1}^d E_{H_0,k}\hat \Delta_k^2 +2\sum_{k=1}^d E_{H_0,k}\hat \Delta_k\r)\\
\nonumber &=& \exp\l(\sum_{k=1}^d E_{0,k}\l( \hat \Delta_k^2
\frac{dP_{H_0,k}}{dP_{0,k}}\r)+o(1) \r)\\ \nonumber
&\le&\exp\l(\frac{p^2}{1-p}\,\sum_{k=1}^d E_{0,k}\l[
L^2(\bY^k)(\ell_{a_d}(Z^k)-1)^2\1_{\{\CA_{SY,k}\bigcap\CA_{SZ,k}\}}
\r]+o(1)\r)
\nonumber\\
 &=&\exp\l(\frac{dp^2A}{1-p}+o(1)\r)
  ,\label{lower11}
\end{eqnarray}
where
$
A=E_{0,1}\l[
L^2(\bY^1)(\ell_{a_d}(Z^1)-1)^2\1_{\{\CA_{SY,1}\bigcap\CA_{SZ,1}\}} \r].
$
Observe that
$$
A\le B+C,\quad B=E_{0,1}\l(
L^2(\bY^1)\ell_{a_d}^2(Z^1)\1_{\{\CA_{SZ,1}\}} \r),\quad C= E_{0,1}\l(
L^2(\bY^1)\1_{\{\CA_{SY,1}\}} \r).
$$
Setting $b_l=a_d \sqrt{l}$ with $l=m \text{ or } m+1,\ T_d=\sqrt{2\log
d}$, we can write
\begin{eqnarray*}
B=\frac{1}{\sqrt{2\pi}}\int_{-\infty}^{T_d}
e^{-b_{m+1}^2+2b_{m+1}t-t^2/2}dt =
 e^{b_{m+1}^2}\Phi(T_d-2b_{m+1}),
\end{eqnarray*}
and analogously,
$$
C=e^{b_{m}^2}\Phi(T_d-2b_{m}).
$$
Recall that we consider $\b\in(3/4,1)$ under assumption
(\ref{lower8}). Thus,
\begin{equation}\label{lower12}
1/2<2\b-1\le (m+1)s^2\le 2(1-\sqrt{1-\b})^2.
\end{equation}
Next, by\nref{ad},
$$
-T_d+2b_{m+1}=\sqrt{2\log d}\,\l(\sqrt{2} \;
s\sqrt{m+1}-1\r)=\sqrt{2\log d}\,\l(\sqrt{2} x_1-1\r).
$$
Thus, for $1/\sqrt{2} < x_1 \le \sqrt{2}$ we have
$$
dp^2B=dp^2e^{b_{m+1}^2}\Phi(T_d-2b_{m+1})\asymp
\frac{d^{-2\b+2\sqrt{2}x_1-x_1^2}}{\sqrt{\log
d}}=\frac{d^{2(1-\b)-(\sqrt{2}-x_1)^2}}{\sqrt{\log d}}.
$$
Here the exponent is $2(1-\b)-\l(\sqrt{2}-x_1\r)^2 \le 0$ in view of
the last inequality in \nref{lower12}. Therefore $ dp^2B=o(1)$ as
$d\to +\infty$. In order to control $dp^2C$ observe that the
function $b\mapsto e^{b^2}\Phi(T-2b)$ is increasing for $b$ large
enough and $T>b$. Therefore $C\le B$ for $d$ large enough and
$dp^2C=o(1)$ as well. Thus $dp^2A=o(1)$ as $d\to +\infty$, and
\nref{lower10a} follows. This completes the proof of
theorem \ref{Lo2}.

%%%%%%%%%%%%%%%%%%%%%%

\subsection{Proof of theorem \ref{Lo2large}}\label{small_l}
Assume w.l.o.g. that $\sigma=1$. By assumptions of the theorem,
$\log m\sim \gamma \log d,\ \g\in (0,1)$ and
\begin{equation}\label{m.b}
\b\in ((1-\g)/2, 1-\g), \quad a=a_d=x\sqrt{\log(d)/m},\quad x=O(1).
\end{equation}
In view of the first two lines of \nref{lower2a}, it suffices to
show that
\begin{equation}\label{m.1}
\sum_{k=1}^d E_{H_{0,k}}\Delta_k^2=dE_{H_{0,1}}\Delta_1^2=o(1).
\end{equation}
Set
$$
\Delta_{{\bf Y}^1}=\frac{pL({\bf Y}^1)}{1-p+pL({\bf Y}^1)}
$$
and observe that
\begin{equation}\label{m.0}
\Delta_{{\bf Y}^1}\le(1-p)^{-1}\min(1,pL({\bf Y}^1)).
\end{equation}
Next, by definition,
\begin{eqnarray*}
L({\bf Y}^1)& = & \exp(-ma^2/2+\sqrt{m}a
SY^1)=d^{-x^2/2+xSY^1/\sqrt{\log d}}, \\
SY^1 & \eqdef & m^{-1/2}\sum_{i=1}^mY^1_i.
\end{eqnarray*}
Take a threshold $H_*=t\sqrt{\log d}$ such that $pL_1(H_*)=1$, i.e.,
$$
-\b-x^2/2+xt=0,\quad t=x/2+\b/x.
$$
Then $pL({\bf Y}^1) <  1$ (respectively, $pL({\bf Y}^1) >  1$) is
equivalent to $SY^1 < H_*$ (respectively, $SY^1 > H_*$).

Since $Z^1, Y^1_i$ are independent and, by the condition
$\limsup_{d\to+\infty}x^*<\phi(\b^*)$, the values $a_d$ are bounded
uniformly in $d$ we have
$E_{H_{0,1}}(\ell_{a_d}(Z^1)-1)^2=e^{a_d^2}-1\le c_0a_d^2 $ where
$c_0$ is a constant. Therefore, using \nref{m.0} we find
\begin{eqnarray*}
E_{H_{0,1}}\Delta_1^2&=&E_{H_{0,1}}(\ell_{a_d}(Z^1)-1)^2
E_P\Delta_{{\bf Y}^1}^2\le c_0 a_d^2E_P\Delta_{{\bf Y}^1}^2\\
&\le& \frac{c_0a^2}{1-p}\l(p^2E_P\l(L^2({\bf Y}^1)\1_{\{pL({\bf
Y}^1)\le 1\}}\r)+P(pL({\bf Y}^1)>1)\r),
\end{eqnarray*}
where  $P=(1-p)P_0+pP_a$, $a=a_d$ for brevity, $P_a$ is the Gaussian
measure with the density $f_a(\cdot)$, cf. (\ref{gden}). Note that
$SY^1\sim\CN(0,1)$ under $P_0$ and $SY^1\sim\CN(0,\sqrt{m}a)$ under
$P_a$. Therefore
\begin{eqnarray*}
E\left(L^2({\bf Y}^1)\1_{\{pL({\bf Y}^1)\le 1\}} \right) & = &
E_{P_0}\left( L^2({\bf Y}^1)\1_{\{pL({\bf Y}^1)\le1\}} \right) + p
E_{P_a}\left(L^2({\bf Y}^1)\1_{\{pL({\bf Y}^1)\le 1\}}
\right),\\
P(pL({\bf Y}^1)>1) & = & (1-p)P_0(SY^1>H_*)+pP_a(SY^1>H_*)\\
 & = & (1-p)c+\lambda(1-p),
\end{eqnarray*}
where
$$
c=\Phi(-H_*)=Ad^{-t^2/2},\quad
\lambda=p\Phi(\sqrt{m}a-H_*)=Ad^{-\b-(t-x)_+^2/2},
$$
and $A$ is a logarithmic factor: $b(\log d)^{-1/2}\le A\le B(\log
d)^{1/2}$ for some positive constants $b,B$. It is easy to see that
$c\le \lambda$ and $c= A\lambda$ as $\sqrt{m}a\le H_*$.

Since $L({\bf Y}^1)=\frac{dP_a}{dP_0}({\bf Y}^1)$ we get
$$
E_{P_0}\left( L^2({\bf Y}^1)\1_{\{pL({\bf Y}^1)\le 1\}} \right) =
\frac{1}{\sqrt{2\pi}}\int_{-\infty}^{H_*} \exp(-ma^2+2xz)dz=
e^{ma^2}\Phi(H_*-2\sqrt{m}a),
$$
and
$$
 pE_{P_a}\left( L^2({\bf Y}^1)\1_{\{pL({\bf Y}^1)\le 1\}}\right) =
pE_{P_0}\left( L^3({\bf Y}^1)\1_{\{pL({\bf Y}^1)\le 1\}} \right) \le
E_{P_0}\left( L^2({\bf Y}^1)\1_{\{pL({\bf Y}^1)\le 1\}}\right).
$$
Therefore
$$
E_{H_{0,1}}\Delta_1^2\le \frac{2a^2(1+o(1))}{1-p}(u+\lambda),
$$
where
$$
u=p^2e^{ma^2}\Phi(H_*-2\sqrt{m}a)=Ad^{-2\b+x^2-(2x-t)_+^2/2}.
$$
 It is easily seen that $u=O(\lambda)$ for $H_*\le \sqrt{m}a$, i.e.,
for
 $t\le x$, which is equivalent to $x^2\ge 2\b$.
Also $\lambda=O(u )$ for $H_*\ge 2\sqrt{m}a$, i.e., for $t\ge 2x$,
which is equivalent to $x^2\le 2\b/3$. If
$\sqrt{m}a<H_*<2\sqrt{m}a$, i.e., if $x<t<2x$, then $u=A\lambda=Ac$;
cf. \cite{IS02a}, pp. 295-296.
 The conditions $x<t<2x$ are equivalent to $
2\b/3\le x^2\le 2\b$. Therefore we get
\begin{equation}\label{m.2}
dE_{H_{0,1}}\Delta_1^2=Ad^{\nu_d},\quad
\nu_d=-\gamma+1+\begin{cases}
-\beta& x^2\ge 2\b,\\
-t^2/2& 2\b/3\le x^2\le 2\b,\\
-2\b+x^2 &0<x^2<2\b/3.
\end{cases}
\end{equation}
Thus the relation \nref{m.1} holds true as
\begin{equation}\label{m.3}
\liminf_{d \to +\infty} \nu_d<0.
\end{equation}
Set
\begin{equation}\label{m.4}
\b^*=\b/(1-\g),\quad x^*=x/\sqrt{1-\g},\quad t^*=x^*/2+\b^*/x^*.
\end{equation}
Then the condition \nref{m.3} is equivalent to $\liminf_{d \to
+\infty} \nu_d^*<0$ where
\begin{equation}\label{m.5}
\nu_d^*=\nu_d/(1-\g)=1+\begin{cases}
-\beta^*&\text{as}\ \ (x^*)^2\ge 2\b^*,\\
-(t^*)^2/2&\text{as}\ \  2\b/3\le (x^*)^2\le 2\b^*,\\
-2\b^*+(x^*)^2 &\text{as}\ \ 0<(x^*)^2<2\b^*/3.
\end{cases}
\end{equation}
The relations \nref{m.5} imply that successful classification is
impossible as $ \limsup_{d \to +\infty} x^*-\phi(\b^*)<0$ where $\phi(\b^*)$ is
defined by \nref{phi} for $\b^*\in (1/2,1)$.

%%%%%%%%%%%%%%%%%

\section{Proof of the upper bounds}
\label{sec6}

In this section we prove theorems \ref{Up1} -- \ref{Up4}. Without
loss of generality, we shall assume throughout that $\sigma=1$. We
shall consider that $s$ is fixed in theorems \ref{Up1} and \ref{Up2}
and that $x$ is fixed in theorems \ref{Up3}, \ref{Up4}. The general
case can be treated in a similar way by passing to subsequences
$s_{d}\to s> 0,\ x_d\to x>0$. Sometimes we shall set for brevity
(and without loss of generality) $c=1$ or $C=1$ where $c$ and $C$
are the constants in the definition of $U_{\b,a_d}$.

\subsection{Proof of theorem \ref{Up1}}\label{U0}

Note first that, for any $\delta>0$, uniformly in $u\in
U_{\beta,a_d}$,
\begin{eqnarray}\label{M1}
&&P^{(u)}(|M_0-h(x)\sqrt{\log d}|>\delta)\to 0,\\
&&P_{H_s}^{(u)}(|M-h(x_s)\sqrt{\log d}|>\delta)\to 0, \quad s=0,1,\label{M3}
\end{eqnarray}
as $d\to +\infty$, where $P^{(u)}$ denotes the distribution of ${\bf
Y}$, the notation $x, x_0, x_1$ is defined in \nref{x} and $
h(t)=\max(\,\sqrt{2}, t+\sqrt{2(1-\b)} \ ). $ Indeed, setting
$T(x)=h(x)\sqrt{\log d}\ge \sqrt{2\log d} $, for any $\delta>0$ we
obtain
\begin{eqnarray*}
P^{(u)}(M_0>T(x)+\delta)&\le& \sum_{k=1}^d
P^{(u)}(SY^k>T(x)+\delta)\le
d\Phi(-T(x)-\delta)\\
&&+\sum_{k=1}^d \varepsilon_k\Phi(a_d \sqrt{m}-T(x)-\delta)\\
&\le &o(1)+Cd^{1-\b}\Phi(-\sqrt{2(1-\b)\log d}-\delta)=o(1)
\end{eqnarray*}
as $d\to +\infty$. Next,
\begin{eqnarray*}
P^{(u)}(M_0<T(x)-\delta) & = & \prod_{k=1}^d(1-P^{(u)}(SY^k\ge
T(x)-\delta)) \\
& \le & \exp\l(-\sum_{k=1}^d P^{(u)}(SY^k \ge
T(x)-\delta)\r),
\end{eqnarray*}
and
\begin{eqnarray*}
\sum_{k=1}^d P^{(u)}(SY^k \ge T(x)-\delta)
%&=&\Big(d-\sum_{k=1}^d
%\varepsilon_k\Big)\Phi(-T(x)+\delta)\\&&+\sum_{k=1}^d
%\varepsilon_k\Phi(a_d \sqrt{m}-T(x)+\delta)\\
\ge (d-Cd^{1-\b})\Phi(-T(x)+\delta)  +
cd^{1-\b}\Phi(a\sqrt{m}-T(x)+\delta).
\end{eqnarray*}
If $h(x)=\sqrt{2}$, then $(d-Cd^{1-\b})\Phi(-T(x)+\delta)$ tends to
$+\infty$ as $d\to +\infty$. If $h(x)>\sqrt{2}$, then
$$
cd^{1-\b}\Phi(a\sqrt{m}-T(x)+\delta)
=cd^{1-\b}\Phi(-\sqrt{2(1-\b)\log d}+\delta)\to +\infty
$$
as $d\to +\infty$. This proves (\ref{M1}). The proof of (\ref{M3})
is analogous.

It follows from \nref{M1}-\nref{M3} that if $x_1\le \phi_2(\b)$
(which is the same as $h(x_1)=\sqrt{2}$, implying
$h(x)=h(x_0)=\sqrt{2}$), then $\Lambda_M<1+\delta$, for any
$\delta>0$ with both $P_{H_0}^{(u)}$ and $P_{H_1}^{(u)}$
probabilities tending to $1$ as $d\to +\infty$. Next, let $x_1>
\phi_2(\b)$. Then $h(x_1)>h(x)\ge h(x_0)$. This yields that
$\Lambda_M<1+\delta$ for any $\delta>0$ with $P_{H_0}^{(u)}$
probability tending to $1$ as $d\to +\infty$. Therefore, (\ref{up11})
and (\ref{up12}) follow. We finally prove (\ref{up13}). Using
(\ref{M1}) and (\ref{M3}) we get that, with $P_{H_1}^{(u)}$
probability tending to $1$ as $d\to +\infty$,
$$
\Lambda_M\ge \frac{h(x_1)\sqrt{\log d} -\delta}{h(x)\sqrt{\log d}
+\delta}>\frac{h(x_1)}{h(x)}\left(1 -\delta\right)-\delta
$$
for any $0<\delta<1$, where the last inequality is satisfied for any
$d\ge 2$.
% such that $h(x)\sqrt{\log d}>1$.
Then (\ref{up13}) holds, since we can always choose a small $c_0$ in
the definition of $\psi^{max}$ and a small $\delta$ such that
$$
\frac{h(x_1)}{h(x)}\left(1 -\delta\right)-\delta> 1+c_0.
$$
Finally, note that all the bounds on the probabilities above are
independent of $u$ and thus the convergence of the probabilities is
uniform in $u\in U_{\b,a_d}$ and in $(a_d, \b)$ such that
${h(x_1)}/{h(x)}$ is bounded away from 1.
%For a fixed $m$, the last
%relation holds true uniformly for $(x_1,\b)$ from a compact set in
%the region of $(x_1,\b): \b\in (0, 1),\ x_1> \phi_2(\b)$.
This completes the proof.

%%%%%%%%%%%%%%%%%%%%%%%%%%%%%

\subsection{Proof of theorem \ref{Up2}}\label{U1}

Fix $u\in U_{\b,a_d}$. We first analyse the expectations and the
variances of the statistics $L_0(t)$ and $L(t)$. Recall that $
\Phi(z)\asymp e^{-z^2/2}/z,\quad z\to+\infty, $ which implies
\begin{equation}\label{A}
\Phi(-tT_d)=A_d d^{-(t_+)^2/2},\quad t_+=\max(t,0).
\end{equation}
Here $A_d$ is a positive factor satisfying $A_d=O(1)$ and $A_d^{-1}=
O(\sqrt{\log d})$ for $t= O(1)$ and $d\to +\infty$. In this proof
and the proof of theorem \ref{Up3} below we assume a weaker
condition: $A_d$ is a quantity depending on $d$ (maybe different on
different occasions) such that
\begin{equation}\label{AA}
 |\log A_d |=o(\log d)
 \end{equation}
 as $d\to +\infty$. Arguing under this
 weaker condition will allow us to get  the proof of theorem
 \ref{Tup3} in parallel with that of theorem
 \ref{Up3}.

 The expectations of $L_0(t)$ and $L(t)$ for
any fixed $u\in U_{\b,a_d}$ and $h\le t \le \sqrt{2} $ satisfy
\begin{eqnarray*}
E^{(u)}L_0(t)&=&d^{1-\b}(\Phi(a_d \sqrt{m}-tT_d)-\Phi(-tT_d))\\
 & = & A_d d^{1-\b}\l(d^{-((t-x)_+)^2/2}-d^{-t^2/2}
\r),\\
E_{H_s}^{(u)}L(t)&=&d^{1-\b}(\Phi((m+s)a_d /\sqrt{m+1}-tT_d)-\Phi(-tT_d))\\
 & = & A_d d^{1-\b}\l(d^{-((t-x_s)_+)^2/2}-d^{-t^2/2} \r),
\end{eqnarray*}
where $s=0,1$ (recall that $h\le t_l \le \sqrt{2} $ for all $t_l$ in
the considered grid). Note that if $x>b$ for some constant $b>0$,
then in view of our assumptions on $h$ we have $ d^{-t^2/2}\le
d^{-(t-x)_+^2/2}/2 $ for $t\ge h $ and all $d$ large enough.
Therefore, for $x,x_s>b$ and all $d$ large enough,
\begin{eqnarray*}
E^{(u)}\Delta_0(t) & = & A_d d^{1/2-\b}d^{t^2/4-(t-x)_+^2/2},\\
E_{H_s}^{(u)}\Delta(t) & = & A_d d^{1/2-\b} d^{t^2/4-(t-x_s)_+^2/2}
,
\end{eqnarray*}
where $s=0,1$. Since the maximum of $t^2/4-((t-x)_+)^2/2$ in $0\le
t\le\sqrt{2}$ is attained either at $t=2x$ when $0<x\le 1/\sqrt{2}$,
or at $t=\sqrt{2}$ when $x> 1/\sqrt{2}$, we have, for $x>b$ and all
$d$ large enough,
\begin{eqnarray}\label{as1}
\max_{h\le t\le \sqrt{2}} E^{(u)}\Delta_0(t)=\begin{cases} A_d d^{1/2-\b+x^2/2}, & x\le 1/\sqrt{2},\\
A_d d^{1-\b-((\sqrt{2}-x)_+)^2/2}, & x>1/\sqrt{2}.
\end{cases}
\end{eqnarray}
Analogously, we have for $s=0,1$, $x_s>b$ and all $d$ large enough:
\begin{eqnarray}\label{as2}
\max_{h\le t\le \sqrt{2}} E_{H_s}^{(u)}\Delta(t)=\begin{cases} A_d d^{1/2-\b+x_s^2/2}, & x_s\le 1/\sqrt{2},\\
A_d d^{1-\b-((\sqrt{2}-x_s)_+)^2/2}, & x_s>1/\sqrt{2}.
\end{cases}
\end{eqnarray}
We shall need the exact asymptotics (\ref{as1}) and (\ref{as2}) only
when $x>b$ and $x_s>b$ for some constants $b>0$. For small $x$ and
$x_s$ it will be enough for our purposes to use the fact that the
right-hand sides of (\ref{as1}) and (\ref{as2}) constitute upper
bounds for the corresponding left-hand sides for all $x,x_s>0$.

We now consider bounds for the corresponding variances:
\begin{eqnarray*}
\Var^{(u)}(L_0(t))&\le&
d\Phi(tT_d)\Phi(-tT_d)+d^{1-\b}\Phi(tT_d-a_d \sqrt{m})\Phi(-tT_d+a_d\sqrt{m})\\
 & \le & A_d ( d^{1-t^2/2}+d^{1-\b-(t-x)^2/2}),\\
\Var_{H_s}^{(u)}(L(t)) & \le &
d\Phi(tT_d)\Phi(-tT_d) \\
& & \mathop{+} d^{1-\b}\Phi(tT_d- (m+s)a_d /\sqrt{m+1})
\Phi(-tT_d+(m+s)a_d /\sqrt{m+1})\\
&\le&A_d ( d^{1-t^2/2}+d^{1-\b-(t-x_s)^2/2}),
\end{eqnarray*}
where $s=0,1$. Since for $x>0$ the maximum of $t^2-(t-x)^2$ in $0\le
t\le\sqrt{2}$ is attained at $t=\sqrt{2}$,
\begin{eqnarray*}
\Var^{(u)}(\Delta_0(t)) &\le& A_d (1+d^{-\b+(t^2-(t-x)^2)/2})\le
A_d (1+d^{1-\b-(\sqrt{2}-x)^2/2})\\
&\le& \begin{cases} A_d, & x\le \phi_2(\b),\\
A_d d^{1-\b-(\sqrt{2}-x)^2/2}, & x>\phi_2(\b),
\end{cases}
\\
\Var_{H_s}^{(u)}(\Delta(t))&\le& A_d (1+d^{-\b+(t^2-(t-x_s)^2)/2})\le
A_d (1+d^{1-\b-(\sqrt{2}-x_s)^2/2})\\
&\le& \begin{cases} A_d, & x_s\le \phi_2(\b),\\
A_d d^{1-\b-(\sqrt{2}-x_s)^2/2}, & x_s>\phi_2(\b),
\end{cases}
\end{eqnarray*}
where $s=0,1$. Take $N_0>0$ such that $N_0^2\asymp T_d\gg N$. By
Chebyshev's inequality, for each $l=1,...,N$, each $u\in U_{\b,a_d}$
and $s=0,1$, with $P_{H_s}^{(u)}$-probability greater than
$1-1/N_0^2$ we have
\begin{eqnarray*}
E_{H_s}^{(u)}\Delta_0(t_l)-N_0\max_{0\le t\le
\sqrt{2}}\sqrt{\Var_{H_s}^{(u)}\Delta_0(t)}\le \Delta_0(t_l) \le
E_{H_s}^{(u)}\Delta_0(t_l)+N_0\max_{0\le t\le
\sqrt{2}}\sqrt{\Var_{H_s}^{(u)}\Delta_0(t)}
\end{eqnarray*}
and these inequalities also valid for $\Delta(\cdot)$ instead of
$\Delta_0(\cdot)$. All these inequalities (with $\Delta(\cdot)$ and
$\Delta_0(\cdot)$) simultaneously hold with probability greater than
$1-2N_0^{-2}N\to 1$ (uniformly in $u\in U_{\b,a_d}$). On this event
of high probability we can evaluate $\Delta_0$ and $\Delta$ by
taking the maxima of the expectations and comparing them with the
maxima of the square root of the variances. Proceeding in this way
and using the bounds obtained above we find:
$$
\Delta_0=\begin{cases} O(A_d), & x\le \phi_1(\b),\ \b \le 3/4 \
\text{or} \ x<\phi_2(\b),\ \b > 3/4,\\
d^{1/2-\b+x^2/2+O(h)}, & 1/\sqrt{2}\ge x>\phi_1(\b),\ \b \le
3/4,\\
d^{1-\b-((\sqrt{2}-x)_+)^2/2+O(h)}, & x>1/\sqrt{2},\ \b \le 3/4 \
\text{or} \ x\ge \phi_2(\b),\ \b > 3/4
\end{cases}
$$
with $P^{(u)}$-probability tending to 1 as $d\to +\infty$, and, for
$s=0,1$ :
$$
\Delta=\begin{cases} O(A_d), & x_s\le \phi_1(\b),\ \b \le 3/4 \
\text{or} \ x_s<\phi_2(\b),\ \b > 3/4,\\
d^{1/2-\b+x_s^2/2+O(h)}, & 1/\sqrt{2}\ge x_s>\phi_1(\b),\ \b \le
3/4,\\
d^{1-\b-((\sqrt{2}-x_s)_+)^2/2+O(h)}, & x_s>1/\sqrt{2},\ \b \le 3/4
\ \text{or} \ x_s\ge \phi_2(\b),\ \b > 3/4
\end{cases}
$$
with $P_{H_s}^{(u)}$-probability tending to 1 as $d\to +\infty$ (the
convergence of all the probabilities is uniform in $u\in
U_{\b,a_d}$). Using these relations we get the following results.
First, $\Lambda^*=o(H)$ with $P_{H_0}^{(u)}$-probability tending to
$1$. Next, $\Lambda^*=o(H)$ with $P_{H_s}^{(u)}$-probability tending
to $1$ (for $s=0,1$) if either $x_1\le \phi_1(\b),\ \b \le 3/4 \
\text{or} \ x_1 <\phi_2(\b),\ \b > 3/4.$ Furthermore, if
$x_0<\sqrt{2}-\tau$ for some small $\tau>0$, and either $x_1>
\phi_1(\b)+\tau,\ \b \le 3/4 \ \text{or} \ x_1>\phi_2(\b) +\tau,\ \b
> 3/4,$ then with $P_{H_1}^{(u)}$-probability tending to $1$ we have
$\Lambda^*\ge d^{c\tau}\gg H$ for some $c>0$. Clearly, the
convergence of all the probabilities here is uniform in $u\in
U_{\b,a_d}$. Thus, the theorem follows.

%%%%%%%%%%%%%%%%%%%%%

\subsection{Proof of theorem \ref{Up3}}\label{U2}

Fix $u\in U_{\b,a_d}$. Let $m=m_d\to +\infty$ such that $\log
m=o(\log d)$.   Observe that $a_d$ cannot be ``too large" in view of
\nref{1*}. Also $a_d$ cannot be ``too small" since $x\eqdef
a_d\sqrt{m/\log d}>\phi(\b) \ge b$ for some $b>0$. In particular,
$a_d d^\delta\to+\infty$, for any $\delta>0$, so that $a_d$
satisfies a condition similar to \nref{AA}:
\begin{equation}\label{upperM.2}
|\log a_d|=o(\log d).
\end{equation}

We first analyse the statistic $\Delta(t)$. Clearly, $
E_{H_0}^{(u)}\Delta(t)=0, $ since $E_{H_0}^{(u)}Z^k=0$ and $Z^k$ and
$SY^k$ are independent. We also have $E_{H_0}^{(u)}(Z^k)^2=1$.
Recalling \nref{A}, we obtain
\begin{eqnarray*}
\Var_{H_0}^{(u)}\Delta(t T_d)&=&\frac{1}{d\Phi(-tT_d)}\l( \sum_{k=1 : \e_k=0}^d\Phi(-tT_d)+ \sum_{k=1 : \e_k=1}^d \Phi(-(t-x)T_d)\r)\\
&=&
1+A_d d^{-\b+t^2/2-((t-x)_+)^2/2},\\
D_0^2(x,\b)&\eqdef&\max_{0<t\le \sqrt{2}}\Var_{H_0}^{(u)}\Delta(t)= 1+A_d d^{1-\b-((\sqrt{2}-x)_+)^2/2}\\
&=&1+A_d \begin{cases}O(1), & x\le \phi_2(\b),\\
d^{1-\b-((\sqrt{2}-x)_+)^2/2}, & x\ge \phi_2(\b).
\end{cases}
\end{eqnarray*}
Here and below $A_d$ is a  factor satisfying \nref{AA}. Next,
$$
E_{H_1}^{(u)}\Delta(t)=\frac{a_d}{\sqrt{d\Phi(-t
T_d)}} \sum_{k=1 : \e_k=1}^d \Phi(-(t-x)T_d)=a_d A_d d^{1/2-\b+t^2/4-((t-x)_+)^2/2},
$$
which yields
\begin{eqnarray*}
 & & E(x,\b)\eqdef\max_{0\le
t\le\sqrt{2}}E_{H_1}^{(u)}\Delta(t)=a_d A_d \begin{cases}
d^{1/2-\b+x^2/2}, & x\le 1/\sqrt{2},\\
d^{1-\b-((\sqrt{2}-x)_+)^2/2}, &  x\ge 1/\sqrt{2}
\end{cases}\\
 & & = a_d A_d \begin{cases}
O(1), & x\le \phi_1(\b),\ \b\le 3/4\ \text{or} \  x\le \phi_2(\b),\ \b\ge 3/4,\\
d^{1/2-\b+x^2/2}, &1/\sqrt{2}\ge x>\phi_1(\b),\ \b\le 3/4,\\
d^{1-\b-((\sqrt{2}-x)_+)^2/2}, &  x>\phi_2(\b),\ \b\ge 3/4\ \text{or}\ x\ge1/\sqrt{2},\ \b\le 3/4.
\end{cases}
\end{eqnarray*}
Analogously,
\begin{eqnarray*}
\Var_{H_1}^{(u)}\Delta(t)&=&\frac{1}{d\Phi(-t T_d)}\Big( \sum_{k=1 : \e_k=0}^d \Phi(-tT_d)+ \sum_{k=1 : \e_k=1}^d \big((\Phi(-(t-x)T_d)\\
&&+a_d^2
\Phi(-(t-x)T_d)\Phi((t-x)T_d)\big)\Big)\\
&=&
1+A_d (1+a_d^2)d^{-\b+t^2/2-((t-x)_+)^2/2},\\
D_1^2(x,\b) & \eqdef & \max_{0<t\le \sqrt{2}}\Var_{H_1}^{(u)}\Delta(t)=
 1+A_d (1+a_d^2)d^{1-\b-((\sqrt{2}-x)_+)^2/2}\\
&=&1+A_d (1+a_d^2)\begin{cases}O(1), & x\le \phi_2(\b),\\
d^{1-\b-((\sqrt{2}-x)_+)^2/2}, & x\ge \phi_2(\b).
\end{cases}
\end{eqnarray*}
%Here the $a^2$ is logarithmic in $d$, cf. \nref{upperM.2a}.
Suppose that, for some small $\tau >0$,
\begin{equation}\label{upperM.1}
\b\in [1/2+\tau,1-\tau],\quad
x\ge \phi(\b)+\tau.
\end{equation}
These relations and the inequality $ 1/2-\b+x^2/2\ge
1-\b-(\sqrt{2}-x)^2/2 $ imply that under \nref{upperM.1} and
\nref{upperM.2} we have, for some $\tau_1>0,\ \tau_2>0$ depending on
$\tau$ in \nref{upperM.1},
$$
D_s(x,\b)\le d^{-\tau_1} E(x,\b), \quad s=0,1 \quad \mbox{and} \quad E(x,\b)\ge
d^{\tau_2}.
$$
Arguing as in the proof of theorem \ref{Up2} above we obtain the
following facts. First,
\begin{equation}\label{upperM.4}
|\Delta| \le A_d D_0(x,\b)
\end{equation}
with $P_{H_0}^{(u)}$-probability tending to $1$ as $d\to +\infty$.
Second, if $x\le\phi(\b)$, then
\begin{equation}\label{upperM.5}
|\Delta|\le A_d D_1(x,\b)
\end{equation}
with $P_{H_1}^{(u)}$-probability tending to $1$ as $d\to +\infty$.
Finally, if \nref{upperM.1} holds, then
\begin{equation}\label{upperM.6} \Delta\ge A_d E(x,\b)
\end{equation}
with $P_{H_1}^{(u)}$-probability tending to $1$ as $d\to +\infty$.
Thus, with $P_{H_0}^{(u)}$-probability tending to $1$, the ratio
$$
\tilde\Lambda(x,\b)=\frac{\Delta}{\sqrt{H+D_0^2(x,\b)}}
$$
is small. The same holds with $P_{H_1}^{(u)}$-probability tending to
$1$ if $x<\phi(\b)$. To finish the proof, we show that these
properties hold also for $\Lambda_\infty^*$ which differs from
$\tilde\Lambda(x,\b)$ only in that we replace $D_0^2(x,\b)$ by
$\Delta_*$ (note that $\tilde\Lambda(x,\b)$ is not a statistic,
since $D_0(x,\b)$ depends on the unknown parameters $x,\b$). The
distribution of $\Delta_*$ is the same under $P_{H_0}^{(u)}$ and
$P_{H_1}^{(u)}$, and depends only on the parameter $u$. We have
\begin{eqnarray*}
E^{(u)}\l( \Delta_* \r) & = & \sum_{k=1 : \e_k=0}^d
\Phi(-\sqrt{2}T_d)+
\sum_{k=1 : \e_k=1}^d \Phi(-(\sqrt{2}-x)T_d) \\
& = & o(1)+A_d d^{1-\b-(\sqrt{2}-x)_+^2/2},\\
\Var^{(u)} \l(\Delta_*\r) & \le &E^{(u)}(\Delta_*).
\end{eqnarray*}
These inequalities yield that $ H+\Delta_*=H+A_d D_0^2(x,\b)$ with
probability tending to $1$, and the statistic $\Lambda^*_\infty$ has
the properties that we have proved for $\tilde\Lambda(x,\b)$.
Finally, note that the convergence of all the probabilities in the
above argument is uniform in $u\in U_{\b,a_d}$. Thus, the theorem
follows.

%%%%%%%%%%%%%%%%%

\subsection{Proof of theorem \ref{Up4}}\label{U4}

For the statistics $L^1(t)$ we have
$$
E_{H_0}^{(u)}L^1(t)=0,\quad E(t) \eqdef
E_{H_1}^{(u)}L^1(t)=a_dd^{1-\b}\Phi(a_d\sqrt{m}-tT_d)=
Aa_dd^{1-\b-(t-x)_+^2/2};
$$
$$
\Var_{H_0}^{(u)}L^1(t)=d(1-p)\Phi(-tT_d)+dp\Phi(a_d\sqrt{m}-tT_d),
$$
$$
\Var_{H_1}^{(u)}L^1(t)=d(1-p)\Phi(-tT_d)+dp(1+a_d^2)\Phi(a_d\sqrt{m}-tT_d),
$$
which yields
$$
\Var_{H_0}^{(u)}L^1(t)\le 2R(t),\quad  \Var_{H_1}^{(u)}L^1(t)\le
3R(t)
$$
with
$$
R(t)=\max\l(d\Phi(-tT_d),\ dp\Phi(a_d\sqrt{m}-tT_d)\r).
$$
Thus, for all $l=1,\dots,N$, with $P^{(u)}_{H_0}$-probability
tending to $1$ the statistics $L^1(t_l)$ belong to the intervals $[-
N\sqrt{2R(t_l)},\ +N\sqrt{2R(t_l)}]$ and with
$P^{(u)}_{H_1}$-probability tending to $1$ they belong to the
intervals $[E(t_l)- N\sqrt{3R(t_l)},\ E(t_l)+N\sqrt{3R(t_l)}]$.

Consider the ratios $\Delta(t_l),\ l=1,...,N$.  First, let
$R(t_l)\le 4N^2$. Then for all $l=1,...,N$, with
$P^{(u)}_{H_s}$-probability tending to $1$ ($s=0,1$), we have the
inequalities
\begin{eqnarray*}
N^2\le N^2+L^0(t_l)\le N^2+2R(t_l)+N\sqrt{R(t_l)} \leq 11N^2,\\
L^1(t_l)\le N\sqrt{2R(t_l)}< 3N^2\quad \text{for}\ s=0,\\
L^1(t_l)\ge E(t_l)-N\sqrt{2R(t_l)}\ge E(t_l)-3N^2 \quad \text{for}\
s=1.
\end{eqnarray*}
Therefore, we get for all $l=1,...,N$ such that $R(t_l)\le 4N^2$,
with $P^{(u)}_{H_s}$-probability tending to $1$,
\begin{eqnarray*}
\Delta(t_l)&<& 4N\quad \text{for}\ s=0,\\
\Delta(t_l)&\ge& E(t_l)/4N-N\ge E(t_l)/(2\sqrt{R(t_l)}\,)-N\quad
\text{for}\ s=1.
\end{eqnarray*}
  Next, let
$R(t_l)> 4N^2$. Then analogously, with $P^{(u)}_{H_s}$-probability
tending to $1$,
\begin{eqnarray*}
N^2+L^0(t_l)\le N^2+2R(t_l)+N\sqrt{R(t_l)}< 3R(t_l),\\
N^2+ L^0(t_l)\ge R(t_l)-N\sqrt{R(t_l)}> R(t_l)/2,\\
L^1(t_l)\le N\sqrt{2R(t_l)}\quad \text{for}\ s=0,\\
L^1(t_l)\ge E(t_l)-N\sqrt{2R(t_l)} \quad \text{for}\ s=1.
\end{eqnarray*}
Hence, we get for all $l=1,...,N$ such that $R(t_l)> 4N^2$, with
$P^{(u)}_{H_s}$-probability tending to $1$,
\begin{eqnarray*}
\Delta(t_l)&<& 4N\quad \text{for}\ s=0,\\
\Delta(t_l)&\ge& E(t_l)/(2\sqrt{R(t_l)}\,)-N\quad \text{for}\ s=1.
\end{eqnarray*}

Thus uniformly over $u\in U_{\b,a_d}$,
$$
E_{H_0}^{(u)}\psi_\infty=P_{H_0}^{(u)}(\Delta>4N)\to 0.
$$

Recalling  \nref{m.b}, \nref{m.4} let us show that under the
condition
\begin{equation}\label{dist}
x^*>\phi(\b^*)
\end{equation}
we have, for some $\eta>0$,
\begin{equation}\label{dist.1}
\max_{1\le l\le N}E(t_l)/\sqrt{R(t_l)}>d^{\eta}.
\end{equation}
This implies that uniformly over $u\in U_{\b, a_d}$,
$$
E_{H_1}^{(u)}(1-\psi_m^*)=P_{H_1}^{(u)}(\Delta\le 4N)\to 0.
$$
In order to verify \nref{dist.1}, let us study the ratio
$E(t)/\sqrt{R(t)}$. We have, with a logarithmic factor $A$,
$$
\frac{E(t)}{\sqrt{R(t)}}=Ad^{s(t)},\quad
s(t)=-\g/2+1/2-\b+\frac{1}{2}\min(t^2/2-(t-x)_+^2,\ \b-(t-x)_+^2/2).
$$
Set $t_0=x/2+\b/x$. Let us  check that
\begin{equation}\label{s*}
s^*\eqdef \max_{0\le t\le\sqrt{2}}s(t)\ge -\frac{\g}{2}+\frac 12+\begin{cases}-\b/2,& \text{if}\ \  x\ge t_0,\\
-t_0^2/4,& \text{if}\ \  x\le t_0\le 2x,\\
-\b+x^2/2,& \text{if}\ \  2x\le t_0.
\end{cases}
\end{equation}
Indeed, the relation $x\ge t_0$ is equivalent to $x^2\ge 2\b $. So,
$s^*\ge s(\sqrt{2\b})=-\g/2+(1-\b)/2$, which implies the first
relation \nref{s*}.

The relation $2x\le t_0$ is equivalent to $x^2\le 2\b/3 $ and if
$2x\le \sqrt{2}$, then $s^*\ge s(2x)=-\g/2+1/2-\b+x^2/2$, which
implies the third relation \nref{s*}. Let us show that the case $2x>
\sqrt{2}, \ x^2\le 2\b/3$ is impossible under \nref{dist}. In fact,
we have
\begin{equation}\label{ineq}
\sqrt{1-\g}\,\phi(\b^*)\ge \phi(\b),\quad 0\le \b\le 1-\g.
\end{equation}
Combining \nref{ineq} and \nref{dist} we find $ x>  \phi(\b) $. It
is easy to see that $x>\phi(\b)$ and $x^2\le 2\b/3 $ only if $\b\le
3/4$. This implies $2x\le \sqrt{2}$.

The relation $x\le t_0\le 2x$ is equivalent to $2\b/3\le x^2\le
2\b$. If $t_0\le \sqrt{2}$, then $s^*\ge
s(\sqrt{2})=-\g/2+1/2-t_0^2/4$, which implies the second relation
\nref{s*}. Let us show that the case $t_0> \sqrt{2},\ 2\b/3\le x^2$
is impossible if $x^*>\phi(\b^*)$. In fact, it is easy to check that
these inequalities are simultaneously satisfied only if
$x\le\phi_2(\b),\ \b\ge 3/4$. However, \nref{dist} and \nref{ineq}
imply
$$
x>  \phi(\b)=\phi_2(\b),\quad \text{for} \ \b\ge 3/4,
$$
a contradiction. By comparing \nref{s*} with \nref{m.5} and
repeating the argument from the end of Subsection \ref{ProofLo2} we
see that \nref{dist} implies $\lim\inf s^*>0$. Since $s(\cdot)$ is a
Lipschitz function, we can replace the maximum over the interval
$[0,\sqrt{2}]$ by the maximum over our grid with step $\delta$,
inducing the error of order $O(\delta)$. This yields \nref{dist.1}.
\endproof

\noindent\sc{Yu.I.~Ingster: St.Petersburg State Electrotechnical University,\\
5, Prof. Popov str., 197376 St.Petersburg, Russia}

\medskip

\noindent\sc{Ch. Pouet: LATP, University of Provence \\
39, rue F. Joliot-Curie, 13453 Marseille cedex 13, France}

\medskip

\noindent\sc{A.B.~Tsybakov: Laboratoire de
Statistique, CREST, Timbre J340 \\
3, av.\ Pierre Larousse, 92240 Malakoff cedex, France \\
and LPMA, University of Paris 6\\
4, Place Jussieu, 75252 Paris cedex 05, France}


\begin{thebibliography}{}
\bibitem{DJ04} Donoho, D. and Jin, J. (2004) Higher criticism for detecting
sparse heterogeneous mixtures. {\it Ann. Statist.} {\bf 32},
962--994.
\bibitem{DJ08a} Donoho, D. and Jin, J. (2008) Higher criticism thresholding:
Optimal feature selection when useful features are rare and weak.
{\it Proc. Nat. Acad. Sci.} {\bf 105}, 14790--14795.
\bibitem{DJ08b} Donoho, D. and Jin, J. (2008) Feature selection by higher
criticism thresholding: Optimal phase diagram. {\it Manucript},
available at arXiv:0812.2263.
\bibitem{HPG} Hall, P., Pittelkow, Y. and Ghosh, M. (2008)
 Theoretical measures of relative
performance of classifiers for high dimensional data with small
sample sizes.{\it J. R. Stat. Soc.} B,{\bf 70}, Part 1, 159--173.
\bibitem{IH} Ibragimov, I. A. and Has'minski, R. Z. (1981)
 {\it Statistical Estimation. Asymptotic Theory.} Springer, New York.
\bibitem{I} Ingster, Yu.I. 1997
 Some problems of hypothesis testing leading to infinitely
divisible distributions. { \it Math. Methods of Stat.} {\bf 6},
47--69.
\bibitem{IS01a} Ingster, Yu. I. and Suslina, I. A. (2001)
Adaptive detection of a signal of growing dimension. I. { \it
Math. Methods of Stat.} {\bf 10}, 395--421.
\bibitem{IS01b} Ingster, Yu. I. and Suslina, I. A. (2001)
Adaptive detection of a signal of growing dimension. II. { \it
Math. Methods of Stat.} {\bf 11}, 37--68.
\bibitem{IS02a} Ingster, Yu. I. and Suslina, I. A. (2002)
{\it Nonparametric Goodness-of-Fit Testing under Gaussian Model}.
Springer Lectures Notes in Statistics.  Vol. {\bf 169}, Springer, New York.
\bibitem{IS02b} Ingster, Yu. I. and Suslina, I .A. (2002)
On a detection of a signal of known shape in multichannel system.
{\it Zapiski Nauchn. Sem. POMI}, {\bf 294},  88--112 (in Russian,
translation in {\it J. Math. Sci.} {\bf 127} (2005), 1, 1723--1736).
\bibitem{jin2009} Jin, J. (2009) Impossibility of successful
classification when useful features are rare and weak. Manuscript.
\bibitem{JW} Jager, L. and Wellner, J. A. (2007)
 Goodness-of-fit tests via phi-divergences.
{\it Ann. Statist.} {\bf 35}, 2018--2053.
\bibitem{Pe} Petrov, V.V. (1995)
{\it Limit Theorems of Probability Theory}. Oxford University Press, Oxford.
\bibitem{Po} Pouet, C. (2008) {\it Quelques contributions \`a la th\'eorie des
tests}. M\'emoire d'habilitation \`a diriger des recherches,
Universit\'e Aix-Marseille 1.
%\bibitem{T} Tsybakov, A.B. (2008) {\it Introduction to Nonparametric
%Estimation}. Springer, New York.
\end{thebibliography}
\end{document}